\documentclass[a4paper,12pt]{scrartcl}
\usepackage[english]{babel}
\usepackage{lmodern}
\usepackage{microtype}
\usepackage{amsmath,amssymb,amsthm,amsfonts,graphicx,etoolbox,scrlayer-scrpage}
\usepackage{mathtools}
\usepackage[noblocks]{authblk}
\usepackage{bbm}
\usepackage[toc,page]{appendix}
\usepackage{xcolor}
\usepackage{enumitem}
\usepackage{pagecolor}

\makeatletter
\patchcmd{\@maketitle}{\huge}{\Large}{}{}
\patchcmd{\abstract}{\quotation}{}{}{}
\AtBeginEnvironment{abstract}{\noindent\ignorespaces}
\AtEndEnvironment{abstract}{\par\mbox{}}

\newcommand{\shorttitle}{\@title}
\makeatother
\setkomafont{author}{\small}
\setkomafont{title}{\rmfamily\bfseries}
\setkomafont{disposition}{\rmfamily\bfseries}

\def\AMS#1{\par\noindent \textbf{AMS subject classification: }#1\par}
\newcommand{\acknowledgements}{\par\mbox{}\par\noindent\textbf{Acknowledgements: }}
\newcommand{\keywords}[1]{\par\noindent\textbf{Keywords: }#1}
\newcommand{\var}{\mathrm{Var}}
\DeclareMathOperator*{\esssup}{ess\,sup}

\theoremstyle{plain}
\newtheorem{theorem}{Theorem}
\newtheorem{lemma}{Lemma}
\newtheorem{definition}{Definition}
\newtheorem{remark}{Remark}

\newtheorem{proposition}{Proposition}
\newtheorem{example}{Example}

\usepackage[
colorlinks,
pdfpagelabels,
pdfstartview = FitH,
bookmarksopen = true,
bookmarksnumbered = true,
linkcolor = blue,
plainpages = false,
hypertexnames = false,
citecolor = red] {hyperref}

\hypersetup{
linktoc=page
}

\renewenvironment{abstract}{\bigskip\noindent\begin{minipage}{\textwidth}\setlength{\parindent}{15pt}\paragraph{Abstract:}}{\end{minipage}}

\begin{document}
\title{On the Impact of Approximation Errors on Extreme Quantile Estimation with
Applications to Functional Data Analysis}

\author[1]{Jaakko Pere\thanks{Corresponding author: jaakko.pere@aalto.fi}}
\author[3]{Benny Avelin}
\author[3]{Valentin Garino}
\author[1]{Pauliina Ilmonen}
\author[2]{Lauri Viitasaari}
\affil[1]{Aalto University School of Science, Finland}
\affil[2]{Aalto University School of Business, Finland}
\affil[3]{Uppsala University Department of Mathematics, Sweden}

\maketitle

\begin{abstract}
	We study the effect of approximation errors in assessing the extreme
	behavior of heavy-tailed random objects. We give conditions for the
	approximation error such that the standard asymptotic results hold for the
	classical Hill estimator and the corresponding extreme quantile estimator.
	As an application, we consider the effect of discretization errors in the
	computation of the $L^p$-norms related to functional data. We approximate
	the norms both with Riemann sums and with Monte Carlo integration. We
	quantify connections between the number of observed functions, the number of
	discretization points, and the regularity of the underlying functions. In
	addition, we derive a new concentration inequality for order statistics.
	This, to the best of our knowledge, is the first Chernoff-type concentration
	inequality for order statistics presented in the literature that provides an
	explicit rate at which the ratio between order statistics and tail quantile
	function converges to one. In our application, the bound is used to provide
	concentration inequalities measuring the distance between the Hill estimator
	based on approximated norms and the one based on the true ones. 
\end{abstract}

\keywords{Approximation error, Concentration, Extreme quantile estimation,
Functional data analysis, Hill estimator}

\smallskip

\AMS{62G32, 60G70, 62R10}

\section{Introduction}
\label{sec:intro}

Heavy-tailed phenomena arise naturally in multiple different application areas
such as in telecommunications, finance and insurance~\cite{resnick2007}. In the
present article, we analyze the effect of approximation errors in assessing the
extreme behavior of heavy-tailed random objects. In particular, we assess the
impact of measurement errors on estimation of the extreme value index $\gamma$
based on the standard Hill estimator and on the estimation of extreme quantiles.
For details on the Hill estimator and the corresponding quantile estimator,
see~\cite{deHaan2007}.

The developed framework is applied in the context of functional data analysis.
Assessment of the effect caused by the approximation error has a prominent role
in functional data analysis. This is due to the fact that one cannot fully
record infinite dimensional observations, and hence true observations have to be
approximated. In addition, often the observed discretizations are noisy as well,
e.g., due to measurement errors. For a review of reconstructing observations in
functional data analysis with interpolation and smoothing methods,
see~\cite{ramsay2005}. 

The literature on the intersection of functional data analysis and extreme value
theory is rather thin. In particular, in related work, we can only
cite~\cite{kim2019}, which examines the consistency of the Hill estimator when
computed on sample scores associated with functional principal component
analysis and assumes heavy-tailed theoretical scores.

In our application, we consider the $L^p$-norms of functional observations $Y_i,
i = 1, \ldots, n$. For instance, with $p = \infty$, this could be applied in
measuring the worst losses of an insurance company over a certain period of
time, in which case extreme losses would mean that capital adequacy is
compromised. As our focus is on quantifying the effect of approximation error,
we assume a priori that the distributions of $L^p$-norms $\|Y_i\|_p$ are
heavy-tailed. We do not consider assessing the heavy-tailedness of the norms
$\|Y\|_p$ -- or more generally some functional $h(Y)$ -- as this separate
problem  has been considered for example in~\cite{hult2005}. 

For the approximation of the norms we use 1) Riemann sums, and 2) Monte Carlo
integration. The first case corresponds, e.g., to the case when the paths
$Y_i(t), t\in[0,1]$, are not fully recorded. The second case corresponds, e.g.,
to the case when the paths $Y_i(t), t\in[0,1]$ are fully recorded, but the
integral cannot be computed explicitly and hence has to be approximated. We
analyze the impact of the approximation error by assessing connections between
the number of observed functions, the number of discretization points, and the
regularity of the underlying functions. This quantifies how refined the
discretization has to be in order to obtain asymptotic results for the Hill
estimator and the corresponding quantile estimator.

In addition, we derive a concentration inequality for order statistics, cf.
Theorem \ref{theorem:concentration}. That is, we derive upper bounds for the
probabilities $\mathbb{P}\left(\left|\frac{X_{n-k, n}}{U(n/k)} - 1\right| > x
\right)$ and $	\mathbb{P}\left(\left|\frac{U(n/k)}{X_{n-k, n}} - 1\right| > x
\right)$, where $X_{n-k,n}$ are the order statistics, $U$ is the tail quantile
function of $X$, and $k=k_n$ is an intermediate sequence. To the best of our
knowledge, these are the first Chernoff-type concentration inequalities studying
the convergence of $\frac{X_{n-k, n}}{U(n/k)}$ towards 1. On related literature,
we mention \cite{aghbalou2024, boucheron2012, fresen2022, xia2019}. In
\cite{aghbalou2024} the authors consider cross-validation on extreme regions, in
\cite{fresen2022} the author studies deviation inequalities for linear
combinations of independent heavy-tailed random variables, and in \cite{xia2019}
the author examines percentile bounds for independent random variables. In
particular, the results presented in \cite{boucheron2012} are closely related to
ours. The authors provide bounds for the moment generating function for the
quantity $X_{n-k,n} - \mathbb{E} X_{n-k,n}$ that translate into (one-sided)
Chernoff-type upper bound for the probability $\mathbb{P}\left(X_{n-k,n} -
\mathbb{E} X_{n-k,n}>x\right)$. 

When considering the Hill estimator based on approximated norms, our Theorem
\ref{theorem:concentration} allows detailed analysis in terms of concentration
inequalities. Namely, Part \ref{part4-continuity} of Theorem
\ref{theorem:continuity} and Part \ref{part4-monte} of Theorem
\ref{theorem:monte} give bounds for the probabilities
$\mathbb{P}\left(|\hat\gamma_n - \tilde\gamma_n| > x\right)$, where
$\hat\gamma_n$ is based on the Hill estimator on approximated values and
$\tilde\gamma_n$ is based on the Hill estimator on the true $L^p$-norms. If one
wishes to compare $\hat\gamma_n$ to the true extreme index $\gamma$, union bound
gives 
\begin{equation*}
	\mathbb{P}\left(|\hat\gamma_n - \gamma| > x\right)
	\leq \mathbb{P}\left(|\hat\gamma_n - \tilde\gamma_n| > \frac{x}{2}\right)
	+ \mathbb{P}\left(|\tilde\gamma_n - \gamma| > \frac{x}{2}\right).
\end{equation*}
Note that concentration inequalities for the estimator $\tilde\gamma_n$ are
considered, under the von Mises' conditions, in~\cite[Theorem
3.3.]{boucheron2015}.

The rest of the article is organized as follows. In Section \ref{sec:general} we
review necessary preliminaries on extreme value theory, and we provide our main
theorems: Theorem \ref{theorem:approx} on the impact of approximation error and
Theorem \ref{theorem:concentration} on the concentration for order statistics.
In Section \ref{sec:fda} we apply the main results in functional data setting.
Section \ref{sec:riemann} considers Riemann sum approximated $L^p$-norms, and
Section \ref{sec:monte} considers Monte Carlo integrated $L^p$-norms. All the
proofs and technical derivations are postponed to the appendix.

\section{General framework}
\label{sec:general}

Let $X_1, \ldots, X_n$ be i.i.d.\ copies of a heavy-tailed random variable $X$
with a cumulative distribution function $F_X$. The tail quantile function $U$
corresponding to the random variable $X$ is defined by
\begin{equation*}
  U(t) = F^\leftarrow\left(1 - \frac{1}{t}\right), \quad t > 1,
\end{equation*}
where $f^\leftarrow$ denotes the left-continuous inverse of a non-decreasing
function $f$. In this article, we utilize the notation $U_X$ to signify that $U$
corresponds to a random variable $X$ when the correspondence is not evident from
the context. Order statistics corresponding to the sample $X_1, \ldots, X_n$ are
denoted by $X_{1,n} \leq \cdots \leq X_{n,n}$.
 
Below we give the definition for regular variation, and three equivalent
definitions for heavy-tailed random variable. 

\begin{definition}[Regular variation]
	A measurable function $f:\mathbb{R}_{>0}\to\mathbb{R}$ that is eventually
	positive is regularly varying with index $\alpha \in \mathbb{R}$, denoted by
	$f\in RV_\alpha$, if
	\begin{equation} \label{eq:regular-variation}
	  	\lim_{t\to\infty}\frac{f(tx)}{f(t)} = x^\alpha, \quad x > 0.
	\end{equation}
	If $f\in RV_0$, we say that $f$ is slowly varying.
\end{definition}

\begin{theorem}[\cite{embrechts1997}, Theorem 3.3.7 and~\cite{deHaan2007},
  	Corollary 1.2.10] \label{theorem:heavy} A random variable $X$ is
  	heavy-tailed with the extreme value index $\gamma > 0$, denoted by $X\in
  	RV_\gamma$, if any of the following equivalent conditions are satisfied:
  	\begin{enumerate}
		\item\label{item:max} Let $X_1, \ldots, X_n$ be i.i.d\ copies of the
		random variable $X$. There exists sequences $a_n > 0$ and $b_n \in
		\mathbb{R}$ such that, as $n\to\infty$,
		\begin{equation*}
	  		\frac{\max_{i\in\{1, \ldots, n\}} X_i - b_n}{a_n}
			\stackrel{d}{\to} G,
		\end{equation*}
		where $G$ is a random variable with the Fr\'echet distribution
		\begin{equation*}
	  		F_G(x) =
	  		\begin{cases}
		  		0, & x\leq 0, \\
		  		\exp\left(-x^{-1/\gamma}\right), & x > 0.
			\end{cases}
		\end{equation*} \label{item:domain}
		\item $1-F_X\in RV_{-1/\gamma}$. \label{item:F-RV}
		\item $U_X\in RV_\gamma$. \label{item:U-RV}
  	\end{enumerate}
\end{theorem}
Note that if $X\in RV_\gamma$, then $x^* = \sup\{x\in\mathbb{R}: F_X(x) < 1\} =
\infty$. Moreover, we can choose $b_n \equiv 0$ and $a_n = U_X(n)$ in
Item~\ref{item:max} of Theorem \ref{theorem:heavy}.

The Hill estimator $\tilde\gamma_n$, introduced in~\cite{hill1975}, is one of
the most well-known estimators for the extreme value index $\gamma > 0$ in the
case of heavy-tailed random variables. Consistency of the Hill estimator was
first established in~\cite{mason1982} under i.i.d.\ heavy-tailed random
variables. The asymptotic normality of the Hill estimator has been studied by
various authors~\cite{haeusler1985, hall1982}. For introduction to extreme value
theory including discussion on different tail index estimators and extreme
quantile estimators, see~\cite{deHaan2007}. The Hill estimator for the tail
index $\gamma$ is defined as 
\begin{equation} \label{eq:hill_q_true}
  	\tilde\gamma_n = \frac{1}{k} \sum_{i=0}^{k - 1}
  	\ln\left(\frac{X_{n - i, n}}{X_{n - k, n}}\right).
 \end{equation}
The corresponding extreme quantile estimator for $x_q =
F^\leftarrow(1 - q) = U(1 / q)$, with $q\in(0,1)$, is defined as 
\begin{equation*}
	\tilde x_q = X_{n - k, n} \left(\frac{k}{nq}\right)^{\tilde\gamma_n}.
\end{equation*}
In asymptotic results for the extreme quantile estimator, it is often assumed
that $q = q_n$ and $q\to 0$ fast as $n\to\infty$. This is necessary since
otherwise quantile $x_q$ would not be extreme for a sufficiently large sample
size $n$. Note that, as is customary in the extreme value theory literature, in
\eqref{eq:hill_q_true} it is implicitly assumed $X_{n-k,n}>0$ since otherwise
$\tilde{\gamma}_n$ is not well-defined. In practice this is the case with
increasing probability when $n\to \infty$.

For the asymptotic normality of the Hill estimator $\tilde\gamma_n$ and the
extreme quantile estimator $\tilde x_q$, conditions presented in Theorem
\ref{theorem:heavy} are not sufficient, see~\cite[Theorem 3.2.5]{deHaan2007}
and~\cite[Theorem 4.3.8]{deHaan2007}. Instead, we pose a second-order regular
variation condition on $X$, denoted by $X\in 2RV_{\gamma, \rho}$, or an extended
second-order regular variation condition, denoted by $X \in 2ERV_{\gamma,
\tilde\rho}$. Here the second-order conditions are connected to $\rho$ and
$\tilde\rho$ that are related to the rate at which the convergence
in~\eqref{eq:regular-variation} occurs, and to the tail index $\gamma$. For
precise definitions, see Definition \ref{def:2rv} and Definition \ref{def:2erv}.

In our setting, we do not observe $X_1, \ldots, X_n$ directly, but instead, we
have approximations $\hat X_1, \ldots, \hat X_n$ that are used in the
estimation. We define the Hill estimator and the extreme quantile estimator
computed on the approximated values as
\begin{equation} \label{eq:hill_q_approx}
  	\hat\gamma_n = \frac{1}{k} \sum_{i=0}^{k - 1}
  	\ln\left(\frac{\hat X_{n - i, n}}{\hat X_{n - k, n}}\right)
  	\quad\text{and}\quad
  	\hat x_q = \hat X_{n - k, n} \left(\frac{k}{nq}\right)^{\hat\gamma_n}.
\end{equation}
Throughout the article, we denote the error by $E_i = \left|\hat X_i -
X_i\right|, i=1,\ldots,n$ and the corresponding order statistics by
$E_{i,n},i=1,\ldots,n$. 
\begin{theorem} \label{theorem:approx} Let $X_1, \ldots, X_n$ be i.i.d.\ copies
  	of $X\in RV_\gamma$, $\hat X_1, \ldots, \hat X_n$ the corresponding
  	approximations, and $\hat\gamma_n$ and $\hat x_q$ given
  	by~\eqref{eq:hill_q_approx}. Assume that, as $n\to \infty$, $k =
  	k_n\to\infty$ and $k/n\to 0$. Then, as $n\to\infty$, we have:
  	\begin{enumerate}
    	\item\label{item:hill-consistency} If $\frac{E_{n,n}}{U(n/k)}
    	\stackrel{\mathbb{P}}{\to} 0$, then
    	$\hat\gamma_n\stackrel{\mathbb{P}}{\to}\gamma$.
 
    	\item \label{item:hill-normality} Assume further that $X\in 2RV_{\gamma,
    	\rho}$ and that $\lim_{n\to\infty} \sqrt{k}A(n/k) = \lambda
    	\in\mathbb{R}$, where $A$ is the auxiliary function from Definition
    	\ref{def:2rv}. If $\sqrt{k}\frac{E_{n,n}}{U(n/k)}
    	\stackrel{\mathbb{P}}{\to} 0$, then
    	\begin{equation*}
      		\sqrt{k}(\hat\gamma_n - \gamma) \stackrel{d}{\to}
	  		N\left(\frac{\lambda}{1 - \rho},\gamma^2\right).
    	\end{equation*}
   	
    	\item \label{item:hill-quantile} Suppose the assumptions of Item 2 hold
    	with $\rho < 0$. Let $q = q_n$, $nq = o(k)$, $\ln(nq) = o(\sqrt{k})$,
    	and denote $d_n = k/(nq)$. If $\sqrt{k}\frac{E_{n,n}}{U(n/k)}
    	\stackrel{\mathbb{P}}{\to} 0$, then
    	\begin{equation*}
      		\frac{\sqrt{k}}{\ln d_n}\left(\frac{\hat x_q}{x_q} - 1\right)
      		\stackrel{d}{\to} N\left(\frac{\lambda}{1 - \rho},\gamma^2\right).
    	\end{equation*}
  	\end{enumerate}
\end{theorem}

\begin{remark} \label{remark:error} In order to control the quantity
	$E_{n,n}/U(n/k)$, note that the tail quantile function can be written as
	$U(x) = L(x) x^\gamma$ for some slowly varying function $L$. It follows from
	\cite[Theorem B.1.6]{deHaan2007} that for any $\varepsilon>0$, as $k =
	k_n\to\infty$ and $k/n\to 0$ when $n\to \infty$, we have
	\begin{equation*}
		\frac{E_{n,n}}{U(n/k)} = o_\mathbb{P}\left(\left(\frac{k}{n}\right)
		^{\gamma-\varepsilon} E_{n,n}\right).
	\end{equation*}
		This gives a simple bound of which rate the maximum error $E_{n,n}$ can
		grow. If in addition we have $X\in 2RV_{\gamma, \rho}$ with $\rho < 0$,
		we have $U(t)/t^\gamma\to c\in(0, \infty)$ (see \cite[page
		49]{deHaan2007}) leading to
	\begin{equation*}
		\frac{E_{n,n}}{U(n/k)} = O_\mathbb{P}\left(\left(\frac{k}{n}\right)
		^{\gamma}E_{n,n}\right).
	\end{equation*}
\end{remark}

We close this section by giving a concentration inequality for $X_{n-k,
n}/U(n/k)$ and its inverse $U(n/k)/X_{n-k,n}$. Note that no second-order
condition is required for the result.
\begin{theorem} \label{theorem:concentration} Let $X_1, \ldots, X_n$ be i.i.d.\
	copies of $X\in RV_\gamma$ and assume that, as $n\to \infty$, $k =
	k_n\to\infty$ and $k/n\to 0$. For $\varepsilon>0$, denote
	\begin{equation} \label{eq:b}
		b_{\gamma, \varepsilon, k}\left(x\right)
		= \exp\left(\left(1 - x^{-1/\gamma}
		- \frac{\ln\left(x\right)}{\gamma} + \varepsilon\right)
		k\right).
	\end{equation}
	\begin{enumerate}
		\item Suppose $0 < x < 1$ and let $\varepsilon>0$ be fixed. Then there
		exists $N =N(\varepsilon,x)$ such that for $n\geq N$ we have 
		\begin{align*}
			\mathbb{P}\left(\left|\frac{X_{n-k, n}}{U(n/k)} - 1\right| > x
			\right) 
			\leq &
			\left(1 + \varepsilon\right)\bigg(
			b_{\gamma, \varepsilon, k}\left(1+x\right) \\
			&+ (1-x)^{-\frac{1}{2\gamma}}
			b_{\gamma, \varepsilon, k}\left(1 - x\right)
			\exp\left((1-x)^{-1/\gamma} \frac{k}{k+1}\right)\bigg)
		\end{align*}
		and
		\begin{equation*}
			\begin{split}
				\mathbb{P}\left(\left|\frac{U(n/k)}{X_{n-k, n}} - 1\right| > x
				\right)
				\leq& \left(1 + \varepsilon\right)\bigg(
				b_{\gamma, \varepsilon, k}\left(\frac{1}{1-x}\right) \\
				&+ (1+x)^{\frac{1}{2\gamma}}
				b_{\gamma, \varepsilon, k}\left(\frac{1}{1+x}\right)
				\exp \left((1+x)^{1/\gamma} \frac{k}{k+1}  \right)
				\bigg) \\
				&+ \mathbb{P}\left(X_{n-k,n} \leq 0\right).
			\end{split}
		\end{equation*}

		\item Suppose $x \geq 1$ and let $\varepsilon>0$ be fixed. Then there
		exists $N =N(\varepsilon,x)$ such that for $n\geq N$ we have 
		\begin{equation*}
			\mathbb{P}\left(\left|\frac{X_{n-k, n}}{U(n/k)} - 1\right| > x
			\right)
			\leq \left(1+\varepsilon\right)
		  	b_{\gamma, \varepsilon, k}\left(1+x\right)
		  	+ \mathbb{P}\left(X_{n-k,n} < 0\right)
		\end{equation*}
		and
		\begin{equation*}
			\begin{split}
				\mathbb{P}\left(\left|\frac{U(n/k)}{X_{n-k, n}} - 1\right| > x
				\right) \leq &
				\left(1 + \varepsilon\right)(1+x)^{\frac{1}{2\gamma}}
				b_{\gamma, \varepsilon, k}\left(\frac{1}{1+x}\right)
				\exp\left((1+x)^{1/\gamma} \frac{k}{k+1}  \right) \\
				&+ \mathbb{P}\left(X_{n-k,n} \leq 0\right).
			\end{split}
		\end{equation*}
	\end{enumerate}
\end{theorem}
Note that for each fixed $x \neq 1$ we have $1-x^{-1/\gamma} -
\frac{\ln(x)}{\gamma} + \varepsilon < 0$ for small enough $\varepsilon$. This
implies exponential decay in $k$ for the terms involving function
$b_{\gamma,\varepsilon,k}$. For example, we have, for $0<x<1$,
\begin{equation} \label{eq:exponential-concentration}
	\mathbb{P}\left(\left|\frac{X_{n-k, n}}{U(n/k)} - 1\right| > x
	\right) \leq C_1 e^{-C_2k}.
\end{equation}
Note that three of the bounds also depend on the probability
$\mathbb{P}\left(X_{n-k,n} \leq 0\right)$. This corresponds to the (small)
probability of the event where the Hill estimator is not even well-defined. In
order to obtain concentration rate in $k_n$, one has to analyze how fast
$\mathbb{P}\left(X_{n-k,n} \leq 0\right)$ decays to zero. However, we can bound
\begin{equation*}
	\mathbb{P}\left(X_{n-k,n} \leq 0\right) \leq \mathbb{P}\left(X_{n-k,n}
	\leq U(n/k)(1-y)\right)
\end{equation*}
with arbitrary $y < 1$ giving exponential decay for $\mathbb{P}\left(X_{n-k,n}
\leq 0\right)$ as well.

\section{Regularly varying $L^p$-norms} 
\label{sec:fda}

In this section we assume that $Y_1, \ldots, Y_n$ are random functions of the
space $L^p([0, 1]^d), d\in\mathbb{N}$ for some given $p\geq 1$. That is,
$Y_1,\ldots,Y_n$ are random functions such that the norm
\begin{equation*}
  	\|y\|_p =
  	\begin{cases}
		\left(\int_{[0, 1]^d}\left|y(t)\right|^p \,\mathrm{d}t\right)
		^\frac{1}{p}, &p\in[1, \infty), \\
		\esssup_{t\in [0, 1]^d} |y(t)|, &p = \infty,
  	\end{cases}
\end{equation*}
where integration and essential supremum are taken with respect to the Lebesque
measure, is almost surely finite. We are interested in the extreme behavior of
the norms $X_i = h(Y_i)$ with $h(y) = \|y\|_p$. We stress that other choices of
the functional $h:L^p([0,1]^d) \mapsto \mathbb{R}$ are possible as long as $h$
has suitable regularity properties, but for the simplicity of our presentation,
we restrict ourselves to the case of the norms.

We consider two cases. In the first case, it is assumed that the functions $Y_i$
are observed only at discrete points $t\in [0,1]^d$ in which case we approximate
the $L^p$-norm using Riemann sums. In this case we restrict to the case $d = 1$
since the performance of approximation with Riemann sums deteriorates quickly as
$d$ increases. In the second case, it is assumed that the functions $Y_i$ are
fully observable, and the norm is computed using Monte Carlo integration.

\subsection{Riemann sum approximation}
\label{sec:riemann}

Let $m\in \mathbb{N} \setminus \{0\}$. We assume that the functions $Y_i\in
L^p([0, 1])$, $1\leq p\leq \infty$ are observed at equidistant points $j/m$,
$j\in\{0, \ldots, m - 1\}$ leading to an approximation 
\begin{equation} \label{eq:approximation}
  	\hat X_i =
  	\begin{cases}
		\left(\frac{1}{m}\sum_{j=0}^{m-1}
		\left|Y_i\left(\frac{j}{m}\right)\right|^p\right)^\frac{1}{p},
		&p\in[1,\infty), \\
		\max_{j\in \{0, \ldots, m-1\}}
		\left|Y_i\left(\frac{j}{m}\right)\right|, &p = \infty,
  	\end{cases}
\end{equation}
of the norm $\Vert Y_i\Vert_p$. Let $\hat\gamma_n$ be the Hill estimator and
$\hat x_q$ be the extreme quantile estimator based on the approximations $\hat
X_i$. Under continuity of the functions $Y_i$, Theorem \ref{theorem:approx} and
Theorem \ref{theorem:concentration} lead to the following. 
 
\begin{theorem} \label{theorem:continuity} Let  $Y_1, \ldots, Y_n$ be i.i.d.\
	copies of $Y\in L^p([0,1]), p\in[1,\infty]$ such that $\Vert Y\Vert_p \in
	RV_\gamma$, let $\hat X_1, \ldots, \hat X_n$ be the corresponding
	approximations given by \eqref{eq:approximation}, and let $\hat\gamma_n$ and
	$\hat x_q$ be given by \eqref{eq:hill_q_approx}. Suppose further that, for
	all $s,t\in [0,1]$, $Y$ satisfies
	\begin{equation} \label{eq:criterion}
		|Y(t) - Y(s)| \leq V \phi(|t - s|) \quad\textnormal{a.s.}
	\end{equation}
	for some random variable $V\in RV_{\gamma'}$ and for some continuous
	non-decreasing function $\phi:\mathbb{R}_{\geq 0}\to \mathbb{R}_{\geq 0}$
	with $\phi(0)=0$. Assume that, as $n\to \infty$, $m=m_n \to \infty$,
	$k=k_n\to \infty$ and $k/n \to 0$. Denote $e_n =
	\phi\left(\frac{1}{m}\right) \frac{U_V(n)}{U_{\|Y\|_p}(n/k)}$. Then, as
	$n\to \infty$, we have: 
 	\begin{enumerate}
		\item \label{part1-continuity} If $e_n \to 0$, then
		$\hat\gamma_n\stackrel{\mathbb{P}}{\to}\gamma$.
	
		\item \label{part2-continuity} Assume further that $\|Y\|_p\in
		2RV_{\gamma, \rho}$ and that $ 	\lim_{n\to\infty} \sqrt{k}A(n/k) =
		\lambda \in\mathbb{R}$, where $A$ is the auxiliary function from
		Definition \ref{def:2rv}. If $ \sqrt{k}e_n \to 0$, then
		\begin{equation*}
			\sqrt{k}(\hat\gamma_n - \gamma) \stackrel{d}{\to}
			N\left(\frac{\lambda}{1 - \rho},\gamma^2\right).
		\end{equation*}
	
		\item \label{part3-continuity} Suppose the assumptions of Item
		\ref{part2-continuity} hold with $\rho < 0$. Let $q = q_n$, $nq = o(k)$,
		$\ln(nq) = o(\sqrt{k})$, and denote $d_n = k/(nq)$. If $ \sqrt{k}e_n \to
		0$, then
		\begin{equation*}
			\frac{\sqrt{k}}{\ln d_n}\left(\frac{\hat x_q}{x_q} - 1\right)
			\stackrel{d}{\to} N\left(\frac{\lambda}{1 - \rho},\gamma^2\right).
		\end{equation*}

		\item \label{part4-continuity} Suppose that $\|Y\|_p > 0$ and $\hat X_i
		> 0$ almost surely for all $i=1,\ldots,n$. Let $V \in 2ERV_{\gamma,
		\rho}$ with $\gamma>0$ and $\rho<0$, and let $\tilde\gamma_n$ be the
		Hill estimator computed on $\|Y_1\|_p, \ldots, \|Y_n\|_p$. Let $x > 0$
		and let $b_{\gamma, \varepsilon, k}$ be given by~\eqref{eq:b}. Choose
		$\varepsilon>0$, $y \geq 1$, and $0 < z < \frac{y}{1 + y}$, and set
		$a_{y, z} = \frac{1}{1+y}+z$. If $e_n\to 0$, then there exists
		$N=N(\varepsilon, y, z)$ such that for $n\geq N$ we have
		\begin{equation} \label{eq:concentration-riemann-2}
			\begin{split}
				\mathbb{P}\left(\left|\hat\gamma_n - \tilde\gamma_n\right|
				> x\right) \leq \left(1 + \varepsilon\right)
				&\left(\left(a_{y,z}\right)^{-1/(2\gamma)}
				b_{\gamma, \varepsilon, k}\left(a_{y,z}\right)
				\exp\left(\left(a_{y,z} \right)^{-1/\gamma} \frac{k}{k+1}\right)
				\right. \\
				&\left. \quad 
				+ \left(\frac{e_n}{\gamma'}\right)^{1/\gamma'}
				\left( \left( \frac{1}{z} \right)^{1/\gamma'}
				+ \left(\frac{2(1+y)}{x} \right)^{1/\gamma'} \right)\right).
			\end{split}
		\end{equation}
	\end{enumerate}
\end{theorem}

The following remarks clarify the assumptions of Theorem
\ref{theorem:continuity} and elaborate the concentration inequality
\eqref{eq:concentration-riemann-2}.

\begin{remark} \label{remark:starting-point} Condition \eqref{eq:criterion} is
	related to the continuity of the function $Y$, and the function $\phi$ is
	called the modulus of continuity. For example, the well-known Kolmogorov
	continuity criteria 
	\begin{equation*}
		\mathbb{E}|Y(t)-Y(s)|^\alpha \leq K|t-s|^{1+\beta}
	\end{equation*}
	would give \eqref{eq:criterion} with $\phi(x) = |x|^\eta$, for any $\eta
	<\beta/\alpha$, and some random variable $V$ having $\alpha$-moments. This
	would correspond to H\"older continuity. Typically, in this case $V$ is more
	heavy-tailed than $\Vert Y\Vert_p$, i.e. $\gamma'\geq \gamma$. More
	generally, if $Y$ is a stochastic process with a deterministic starting
	point $Y(0)$ such that \eqref{eq:criterion} holds, then by the Minkowski
	inequality we obtain
	\begin{equation*}
		\|Y\|_p \leq V \|\phi\|_p + |Y(0)|
	\end{equation*}
	from which $\gamma'\geq \gamma$ follows. However, it can also happen that
	$\gamma'<\gamma$. Let, for example, $Y(t) = \mathcal{R} +
	\left|B(t)\right|$, where $\mathcal{R}$ is a non-negative heavy-tailed
	random variable and $B$ is a Brownian motion. Then $\Vert Z\Vert_p \geq
	\mathcal{R}$, and it follows that $\Vert Z \Vert_p$ is at least as
	heavy-tailed as $\mathcal{R}$, while for the increments we have 
	\begin{equation*}
		|Y(t) - Y(s)| \leq |B(t)-B(s)| \leq V|t-s|^{\eta}
	\end{equation*}
	with any $0 < \eta < \frac{1}{2}$. It can be shown that here $V$ has finite
	exponential moments, and hence is lighter-tailed than $\mathcal{R}$. 
\end{remark}

\begin{remark} \label{remark:continuity} The error $e_n =
	\phi\left(\frac{1}{m}\right) \frac{U_V(n)}{U_{\|Y\|_p}(n/k)}$ determines how
	many discretization points $m$ are needed for the approximation error
	staying sufficiently small. In particular, more points are needed for more
	irregular functions $Y$ (giving worse modulus of continuity $\phi$). Note
	also that, similarly as in Remark \ref{remark:error}, we obtain, for any
	$\varepsilon > 0$, 
	\begin{equation*}
		\phi\left(\frac{1}{m}\right) \frac{U_V(n)}{U_{\|Y\|_p}(n/k)}
		= o\left(\phi\left(\frac{1}{m}\right) k^{\gamma - \varepsilon}
		n^{\gamma' - \gamma + 2\varepsilon}\right).
	\end{equation*}
	This also shows that if $V$ is very heavy-tailed compared to $\|Y\|_p$,
	a finer discretization is required so that $\phi\left(\frac{1}{m}\right)$
	compensates the growth of $n^{\gamma'-\gamma}$. Similarly,
	$\phi\left(\frac{1}{m}\right)$ needs to compensate the growth of $k^\gamma$
	that depends on the choice of intermediate sequence $k= k_n$ and the index
	$\gamma>0$. 
\end{remark}

\begin{remark} \label{remark:choice} The upper bound in Concentration
	inequality~\eqref{eq:concentration-riemann-2} consists of two terms. The
	first term comes from the concentration of the order statistic in the event
	that the approximation error is small, while the second term corresponds to
	the probability that the approximation error is large. Moreover, the choice
	of $y$ and $z$ in~\eqref{eq:concentration-riemann-2} is made to
	balance the terms, and allows for some flexibility. Specifically, if $e_n$
	is very small, we can allow for a small $z$ and a large $y$, which
	will make the first term in~\eqref{eq:concentration-riemann-2} small at the
	price of making the second term large. Without knowing anything about the
	specific size of $e_n$, we can choose as follows. Let $0 < r_0 < 1$ be such
	that $1-r_0^{-1/\gamma}-\ln(r_0)/\gamma = -1$. Then, for $0 < r \leq r_0$,
	\begin{equation*}
		b_{\gamma, \varepsilon, k}(r) \leq e^{-k/2}
	\end{equation*}
	for any $0 < \varepsilon < 1/4$. By choosing $\frac{1}{1+y} = r_0/2$ and
	$z = r_0/2$, we have  $0 < z < \frac{y}{1 + y}$, $1 < y$, and
	$\frac{1}{1+y} + z = r_0$. With these choices
	\begin{equation*}
		\mathbb{P}\left(\left|\hat\gamma_n - \tilde\gamma_n\right| > x \right) 
		\leq 
		2 \left(c_0 e^{-k/2}+c_1 e_n^{1/\gamma'}
		\left(1 + \frac{1}{x^{1/\gamma'}} \right)\right)
	\end{equation*}
	for constants $c_0 = c_0(\gamma)>0$ and $c_1= c_1(\gamma,\gamma') > 0$. This
	gives exponential decay in $k$ on the first term
	(cf.~\eqref{eq:exponential-concentration}) and polynomial decay in the error
	$e_n$. Finally, for $0 < y < 1$ we obtain, by using Proposition
	\ref{prop:min-concentration}, that 
	\begin{equation*}
		\begin{split}
			\mathbb{P}\left(\left|\hat\gamma_n - \tilde\gamma_n\right| > x
			\right) \leq \left(1 + \varepsilon\right)
			&\left(b_{\gamma, \varepsilon, k}\left(\frac{1}{1-y}
			- z\right) \right.\\
			&\quad + \left(a_{y,z}\right)^{-1/(2\gamma)}
			b_{\gamma, \varepsilon, k}\left(a_{y,z}\right)
			\exp\left(\left(a_{y,z} \right)^{-1/\gamma} \frac{k}{k+1}\right) \\
			&\left. \quad
			+ \left(\frac{e_n}{\gamma'}\right)^{1/\gamma'}
			\left( \left( \frac{1}{z} \right)^{1/\gamma'}
			+ \left(\frac{2(1+y)}{x} \right)^{1/\gamma'} \right)\right).
		\end{split}
	\end{equation*}
	This is a slight improvement over~\eqref{eq:concentration-riemann-2} in the
	case that $e_n$ is large, as it allows to choose a small $y$.
\end{remark}

A natural class of suitable stochastic processes satisfying
Assumption~\eqref{eq:criterion} of Theorem \ref{theorem:continuity} involves
series $Y(t) = \sum_{j=1}^\infty\varphi_j\mathcal{R}_jZ_j(t)$, where $\mathcal
R_j$ are i.i.d.\ heavy-tailed random variables, $Z_j$ are almost surely
continuous processes having lighter tails than $\mathcal R_j$, and $\varphi_j$
are real constants such that the aforementioned series converges. 

\begin{example} \label{example:riemann} Let $\mathcal{R}_1, \mathcal{R}_2,
	\ldots$ be i.i.d.\ copies of a nonnegative $\mathcal{R}\in RV_{\gamma'}$ and
	let $Z_1, Z_2, \ldots$ be i.i.d. copies of $Z$ satisfying 
    \begin{equation} \label{eq:moments}
    	\mathbb{E}\left(\left| Z(t) - Z(s)\right|^{\kappa}\right)
    	\leq K\left| t - s \right|^{1 + \beta} \quad\forall \ t,s\in[0,1],
    \end{equation}
    where $\beta > 0$, $K > 0$, and $\kappa> 1/\gamma'$. Such processes $Z$
    include, e.g., H\"older continuous Gaussian processes or, more generally,
    hypercontractive H\"older continuous processes \cite{nummi2024}. Assume also
    that $\{\mathcal{R}_j\}_{j\in\mathbb{N}}$ and $\{Z_j\}_{j\in\mathbb{N}}$ are
    mutually independent. Then we can construct $Y$ satisfying
    \eqref{eq:criterion} by setting
    \begin{equation*}
    	Y(t) = \sum_{j = 1}^\infty \varphi_j\mathcal{R}_j Z_j(t),
    	\quad t\in[0, 1],
    \end{equation*}
    where $\{\varphi_j\}_{j\in\mathbb{N}}$ is a (non-trivial) deterministic
    sequence with $\sum_{j = 1}^\infty \left|\varphi_j\right|^\lambda < \infty$
    for some $0 < \lambda < \min\left\{1/\gamma', 1\right\}$. Now if $Z$ is not
    degenerate in the sense that $Z(t) \neq Z(s)$ almost surely for some
    $t,s\in[0,1]$, Then the Hölder constant of $Z$ is almost surely positive,
    and it follows that 
    \begin{equation} \label{eq:criterion_example}
    	|Y(t) - Y(s)| \leq V |t-s|^\eta \quad\textrm{for all}\quad
    	\eta\in (0, \beta / \kappa),
    \end{equation}
    where $V\in RV_{\gamma'}$ (the details on this fact are given in Appendix
    \ref{appendix:sec:examples}). Note that in this case we have $\phi(\delta) =
    |\delta|^{\beta/\kappa - \varepsilon'}$ for any $0 < \varepsilon' <
    \beta/\kappa$. Hence, by choosing $k \sim n^\lambda$ for $\lambda\in(0,1)$
    we, in view of Remark \ref{remark:continuity}, obtain consistency for the
    Hill estimator $\hat\gamma_n$ provided that
    $\left(\frac{1}{m}\right)^{\frac{\beta}{\kappa} - \varepsilon'}n^{\gamma' -
    (1-\lambda)\gamma + \varepsilon'}=O(1)$. This can be obtained by fine-tuning
    $m$ and $\lambda$. 
\end{example}

\subsection{Monte Carlo integration}
\label{sec:monte}

Let $m\in \mathbb{N} \setminus \{0\}$ and $Y_i\in L^p([0, 1]^d)$ for $d\in
\mathbb{N}$ and $1\leq p < \infty$. The norms $X_i = \|Y_i\|_p$ are computed
applying Monte Carlo approximation given by 
\begin{equation} \label{eq:monte}
	\hat X_i = \left(\frac{1}{m}\sum_{j=1}^{m}
	\left|Y_i(T_j)\right|^p\right)^{1/p},
\end{equation}
where $T_1, \ldots, T_m$ are i.i.d.\ copies of a random variable $T \sim
U([0,1]^d)$. Instead of continuity, we pose additional moment assumptions.
Theorem \ref{theorem:approx} and Theorem \ref{theorem:concentration} then lead
to the following. 
\begin{theorem} \label{theorem:monte} Let $Y_1, \ldots, Y_n$ be i.i.d.\ copies
	of $Y\in L^p([0,1]^d), d\in \mathbb{N}$, and $p\in[1,\infty)$ such that
	$\Vert Y\Vert_p \in RV_\gamma$, let $\hat X_1, \ldots, \hat X_n$ be the
	corresponding approximations given by \eqref{eq:monte}, and let
	$\hat\gamma_n$ and $\hat x_q$ be given by \eqref{eq:hill_q_approx}. Suppose
	further that
	\begin{equation} \label{eq:criterion-monte}
		\var\left(\left|Y(T)\right|^p | Y\right) \leq V^{2p}
	\end{equation}
	almost surely for some non-negative $V\in RV_{\gamma'}$ with $\gamma' <
	\frac{1}{2p}$. Assume that, as $n\to \infty$, $m=m_n \to \infty$, $k=k_n\to
	\infty$ and $k/n \to 0$. Denote $e_n' =
	\left(\frac{n}{m}\right)^{\frac{1}{2p}} \frac{U_V(n)}{U_{\|Y\|_p}(n/k)}$.
	Then, as $n\to \infty$, we have: 
	\begin{enumerate}
		\item\label{part1-monte} If $e_n' \to 0$, then
		$\hat\gamma_n\stackrel{\mathbb{P}}{\to}\gamma$.
	
		\item\label{part2-monte} Assume further that $\|Y\|_p\in 2RV_{\gamma,
		\rho}$ and that $\lim_{n\to\infty} \sqrt{k}A(n/k) = \lambda
		\in\mathbb{R}$, where $A$ is the auxiliary function from Definition
		\ref{def:2rv}. If $ \sqrt{k}e_n' \to 0$, then
		\begin{equation*}
			\sqrt{k}(\hat\gamma_n - \gamma) \stackrel{d}{\to}
			N\left(\frac{\lambda}{1 - \rho},\gamma^2\right).
		\end{equation*}

		\item\label{part3-monte} Suppose the assumptions of Item 2 hold with
		$\rho < 0$. Let $q = q_n$, $nq = o(k)$, $\ln(nq) = o(\sqrt{k})$, and
		denote $d_n = k/(nq)$. If $ \sqrt{k}e_n' \to 0$, then
		\begin{equation*}
			\frac{\sqrt{k}}{\ln d_n}\left(\frac{\hat x_q}{x_q} - 1\right)
			\stackrel{d}{\to} N\left(\frac{\lambda}{1 - \rho},\gamma^2\right).
		\end{equation*}
	
		\item\label{part4-monte} Suppose that $\|Y\|_p > 0$ and $\hat X_i > 0$
		almost surely for all $i=1,\ldots,n$. Let $\tilde\gamma_n$ be the Hill
		estimator computed on $\|Y_1\|_p, \ldots, \|Y_n\|_p$. Let $x > 0$ and
		let $b_{\gamma, \varepsilon, k}$ be given by~\eqref{eq:b}. Choose
		$\varepsilon>0$, $y\geq 1$, and $0 < z < \frac{y}{1 + y}$, and
		set $a_{y,z} = \frac{1}{1+y}+z$. Then there exists
		$N=N(\varepsilon,y,z)$ such that for $n\geq N$ we have
		\begin{equation} \label{eq:concentration-monte-2}
			\begin{split}
				\mathbb{P}\left(\left|\hat\gamma_n - \tilde\gamma_n\right| > x
				\right) \leq \left(1 + \varepsilon\right)
				&\left(a_{y,z}^{-1/(2\gamma)}
				b_{\gamma, \varepsilon, k}\left(a_{y,z}\right)
				\exp\left( a_{y,z}^{-1/\gamma} \left( \frac{k}{k+1}\right)
				\right)\right. \\
				&\left. \quad + \mathbb{E}\left(G^{2p}\right) (e_n')^{2p}
				\left(\left(\frac{1}{z}\right)^{2p}
				+ \left(\frac{2\left(1 + y\right)}{x}\right)^{2p}\right)
				\right),
			\end{split}
		\end{equation}
		where $G$ is a Fr\'echet distributed random variable.
	\end{enumerate}
\end{theorem}
As in Section \ref{sec:riemann}, we close this section by clarifying the
assumptions and providing simple examples.

\begin{remark} \label{remark:monte} Condition \eqref{eq:criterion-monte} can be
	equivalently given as
	\begin{equation*}
		\left(\left\|Y\right\|_{2p}^{2p} - \left\|Y\right\|_p^{2p}\right)
		^{\frac{1}{2p}} \leq V
		\quad a.s.
	\end{equation*}
	Indeed, we have $\|Y\|_p =  \left(\mathbb{E}\left(|Y(T)|^p |
	Y\right)\right)^{1/p}$ with $T \sim U([0,1]^d)$. This leads to 
	\begin{equation*}
		\var\left(\left|Y(T)\right|^p | Y\right)
		= \mathbb{E}\left(\left|Y(T)\right|^{2p} | Y\right)
		- \left(\mathbb{E}\left(\left|Y(T)\right|^p | Y\right)\right)^2
		= \left\|Y\right\|_{2p}^{2p} - \left\|Y\right\|_p^{2p}.
	\end{equation*}
	Hence one can choose, for example, $V = \left\|Y\right\|_{2p}$ if
	$\left\|Y\right\|_{2p}\in RV_{\gamma'}$ in which case $\gamma'\geq \gamma$.
	On the other hand, we also have $\left(\left\|Y\right\|_{2p}^{2p} -
	\left\|Y\right\|_p^{2p}\right)^{\frac{1}{2p}} \leq
	2^{\frac{1}{2p}}\left(\left\|Y\right\|_{2p} - \left\|Y\right\|_{p}\right)$
	allowing the choice $V = 2^{\frac{1}{2p}}\left(\left\|Y\right\|_{2p} -
	\left\|Y\right\|_{p}\right)$. Now, as in the case of Riemann sum
	approximations (cf. Remark \ref{remark:starting-point}), it may happen that
	$\gamma'<\gamma$. As an example, consider 
	\begin{equation*}
		Y(t) = B\mathcal{R}_1 + (1 - B)\mathcal{R}_2f(t),
	\end{equation*}
	where $B$, $\mathcal{R}_1$ and $\mathcal{R}_2$ are independent random
	variables such that $\mathbb{P}\left(B = 1\right) = 0.5 = \mathbb{P}\left(B
	= 0\right)$, $\mathbb{P}\left(\mathcal{R}_1 > x\right) = x^{-1/\gamma}$ and
	$\mathbb{P}\left(\mathcal{R}_2 > x\right) = x^{-1/\gamma'}$, $x \geq 1$ for
	$\gamma > \gamma' > 0$, and $f\in L^{2p}([0, 1]^d)$ is a deterministic
	function. It follows that $V=2^{\frac{1}{2p}}\left(\left\|Y\right\|_{2p} -
	\left\|Y\right\|_{p}\right) \in RV_{\gamma'}$ while $\Vert Y\Vert_p \in
	RV_\gamma$.
\end{remark}

\begin{remark} \label{remark:monte2} The error $e_n'=
	\left(\frac{n}{m}\right)^{\frac{1}{2p}} \frac{U_V(n)}{U_{\|Y\|_p}(n/k)}$
	determines how many Monte Carlo approximation points $m$ are needed for the
	approximation error staying sufficiently small. Indeed, similarly as in
	Remark \ref{remark:continuity}, the error satisfies 
	\begin{equation*}
		e_n' = o\left(m^{-1/2p}k^{\gamma-\varepsilon}
		n^{1/2p+\gamma' - \gamma + 2\varepsilon}\right).
	\end{equation*}
\end{remark}

\begin{remark}
	As in Remark \ref{remark:choice}, we can simplify the concentration
	inequalities of Theorem \ref{theorem:monte} by a specific choice of
	$a_{y,z}$. This leads to 
	\begin{equation*}
		\mathbb{P}\left(\left|\hat\gamma_n - \tilde\gamma_n\right| > x
		\right) \leq 2 \left(c_0 e^{-k/2} + c_1 \mathbb{E} (G^{2p}) (e_n')^{2p}
		(1+1/x^{2p}) \right)
	\end{equation*}
	for constants $c_0(\gamma), c_1(\gamma,p) > 0$. Furthermore,
	\eqref{eq:concentration-monte-2} can be slightly improved for large $e_n$ by
	considering the case $0<y<1$.
\end{remark}

\begin{remark} \label{remark:positive-approx} In Part 4 of both Theorem
	\ref{theorem:continuity} and Theorem \ref{theorem:monte} we assume that
	both, the norms $\|Y_i\|_p$ and the approximations $\hat{X}_i$, are almost
	surely positive. Due to continuity, this is always true eventually (by
	taking fine enough discretization) in the case of Riemann sum
	approximations. In the case of Monte Carlo integration this is not
	necessarily true, and one can construct examples where the approximations
	are zero with positive probability. This is not a problem as, naturally,
	such approximations are not used in practice.
\end{remark}

\begin{example} \label{example:monte} Let $\mathcal{R} \in RV_{\gamma'}$ with
	$\gamma'<\frac{1}{2p}$ and let $Z$ be a process satisfying
	$\mathbb{E}\left(\left\|Z\right\|_{2p}^{\kappa}\right) < \infty$ for some
	$\kappa>\frac{1}{\gamma'}$. Such processes $Z$ include, e.g.
	hypercontractive continuous processes mentioned in Example
	\ref{example:riemann} and, as non-continuous examples, L\'evy processes with
	suitable moment conditions. Assume that $\mathcal{R}$ and $Z$ are
	independent. Then we can construct $Y$ satisfying \eqref{eq:criterion-monte}
	by setting $Y(t) = \mathcal{R}Z(t)$. Now if $Z$ is non-degenerate in the
	sense that $\mathbb{P}\left(Z \equiv 0\right) < 1$, then the moment
	$\mathbb{E}\|Z\|_{2p}^{1/\gamma'}$ is positive, and it follows that
	$\|Y\|_{2p}\in RV_{\gamma'}$ (the details on this are given in Appendix
	\ref{appendix:sec:examples}).
\end{example}

\bigskip
\acknowledgements{Jaakko Pere gratefully acknowledges support from the Vilho,
Yrj\"o and Kalle V\"ais\"al\"a Foundation and the hospitality of the Uppsala
University Department of Mathematics during his research visit. Pauliina Ilmonen
gratefully acknowledges support from the Academy of Finland via the Centre of
Excellence in Randomness and Structures, decision number 346308. Benny Avelin
gratefully acknowledges support from the Swedish Research Council, grant number
2019--04098.}

\bibliographystyle{abbrv}
\bibliography{sources}

\begin{appendices}
	\counterwithin{theorem}{section}
	\counterwithin{definition}{section}
	\counterwithin{lemma}{section}
	\counterwithin{remark}{section}
	\counterwithin{example}{section}
	\counterwithin{corollary}{section}
	\counterwithin{proposition}{section}

	In Appendix \ref{appendix:sec:second_order} we present tail regularity
	conditions that are used in the article. The proofs of the theorems are
	presented in Appendices \ref{appendix:sec:general} and
	\ref{appendix:sec:fda}. Details related to Examples \ref{example:riemann}
	and \ref{example:monte} are given in Appendix \ref{appendix:sec:examples}.

	\section{Second-order conditions}
	\label{appendix:sec:second_order}

	Below we give two different second-order conditions that are commonly used
	in asymptotic results in extreme value theory. The first one is used in the
	heavy-tailed case and the second one in a more general framework. For
	details about the relationship between the second-order conditions,
	see~\cite{alves2007, neves2009}.

	\begin{definition}[Second-order regular variation]
		\label{def:2rv}
		A random variable $X$ satisfies the second-order regular variation
		condition with parameters $\gamma > 0$ and $\rho \leq 0$, denoted by
		$X\in 2RV_{\gamma, \rho}$, if the corresponding tail quantile function
		$U$ satisfies
		\begin{equation*}
			\lim_{t\to\infty}\frac{\frac{U(tx)}{U(t)} - x^\gamma}{A(t)}
			= x^\gamma\frac{x^\rho - 1}{\rho}
		\end{equation*}
		 for all $x>0$ and for some positive or negative function $A$ with
		$\lim_{t\to\infty}A(t) = 0$. For $\rho = 0$, the right-hand side is
		interpreted as $x^\gamma\ln x$.
	\end{definition}

	\begin{definition}[Second-order extended regular variation]
	\label{def:2erv}
		A random variable $X$ satisfies the second-order extended regular
		variation condition with parameters $\gamma\in\mathbb{R}$ and
		$\tilde\rho \leq 0$, denoted by $X\in 2ERV_{\gamma, \tilde\rho}$, if
		the corresponding tail quantile function $U$ satisfies
		\begin{equation*}
			\lim_{t\to\infty}\frac{\frac{U(tx) - U(t)}{a(t)}
			- \frac{x^\gamma - 1}{\gamma}}{\tilde A(t)}
			= \frac{1}{\tilde\rho}\left(\frac{x^{\gamma + \tilde\rho} - 1}
			{\gamma + \tilde\rho} - \frac{x^\gamma - 1}{\gamma}\right)
			=: H_{\gamma, \tilde\rho}(x)
		\end{equation*}
		for all $x > 0$, for some positive function $a$, and for some positive
		or negative function $\tilde{A}$ with $\lim_{t\to\infty}\tilde{A}(t) =
		0$. For $\gamma = 0$ or $\tilde\rho = 0$, the right-hand side is
		interpreted as the limit of $H_{\gamma, \tilde\rho}(x)$, as $\gamma\to
		0$ or $\tilde\rho\to 0$.
	\end{definition}

	\section{Proofs of Theorem \ref{theorem:approx} and Theorem
	\ref{theorem:concentration}}
	\label{appendix:sec:general}
	
	\begin{proof}[Proof of Theorem \ref{theorem:approx}]
	    Suppose that the assumptions of Theorem \ref{theorem:approx} hold. 
	    
	    We begin by analyzing the error in order statistics. Let $\ell\in\{n -
		k, \ldots, n\}$. For all $i\in \{1, \ldots, n\}$ we have $|\hat X_i -
		X_i| \leq E_{n,n}$, and hence $X_i - E_{n,n} \leq \hat X_i \leq X_i +
		E_{n,n}$. Denote sets of indices by $I_\ell = \left\{ i \in \left\{1,
		\ldots, n\right\} : X_i \leq X_{\ell,n} \right\}$ and $J_\ell = \left\{
		j \in \left\{1, \ldots, n\right\} : X_j \geq X_{\ell,n} \right\}.$ Now,
		we have 
		\begin{equation*}
			\begin{cases}
				\hat X_i \leq X_{\ell, n} + E_{n,n} & \forall i\in I_\ell \\
				\hat X_j \geq X_{\ell, n} - E_{n,n} & \forall j\in J_\ell
			\end{cases}.
		\end{equation*}
		Since $\# I_\ell = \ell$ and $\# J_\ell = n - \ell + 1$, it follows that
		\begin{equation*}
			\begin{cases}
				\hat X_{\ell, n} \leq X_{\ell, n} + E_{n,n} \\
				\hat X_{\ell, n} \geq X_{\ell, n} - E_{n,n}
			\end{cases}
		\end{equation*}
		and thus
		\begin{equation} \label{eq:order_stat_diff_ineq}
			\left|\hat X_{\ell, n} - X_{\ell, n}\right| \leq E_{n,n}.
		\end{equation}
		If $\hat X_{\ell, n} \geq X_{\ell, n} > 0$, inequality $\ln z \leq z -
		1$ gives 
		\begin{equation*} 
			\left|\ln\left(\frac{\hat X_{\ell, n}}{X_{\ell, n}}\right)\right|
			= \ln\left(\frac{\hat X_{\ell, n}}{X_{\ell, n}}\right)
			\leq \frac{\left|\hat X_{\ell, n} - X_{\ell, n}\right|}{X_{\ell, n}}
			\leq \frac{E_{n,n}}{X_{n-k, n}}.
		\end{equation*}
		Handling the case  $X_{\ell, n} \geq \hat X_{\ell, n}$ similarly leads
		to 
		\begin{equation} \label{eq:ln_ineq}
			\left|\ln\left(\frac{\hat X_{\ell, n}}{X_{\ell, n}}\right)\right|
			\leq \frac{E_{n,n}}{\min\left(\hat X_{n-k, n},
			X_{n-k, n}\right)}.
		\end{equation}
		
		Next, we turn our focus to the Hill estimator. Recall that
		$\left|\frac{X_{n-k,n}}{U(n/k)} - 1\right|\stackrel{\mathbb{P}}{\to}0$
		as $n\to\infty$. Since 
		\begin{equation*}
			\left|\frac{\hat X_{n-k,n}}{U(n/k)} - 1\right|
			\leq \frac{\left|\hat X_{n-k,n} - X_{n-k,n}\right|}{U(n/k)}
			+ \left|\frac{X_{n-k,n}}{U(n/k)} - 1\right|
			\leq \frac{E_{n,n}}{U(n/k)} +
			\left|\frac{X_{n-k,n}}{U(n/k)} - 1\right|
		\end{equation*}
		and since $E_{n,n}/U(n/k)\stackrel{\mathbb{P}}{\to}0$ by assumption, we
		have $\hat X_{n-k,n}/U(n/k)\stackrel{\mathbb{P}}{\to}1$. Note that then,
		as $n\to \infty$, $\max\left(\mathbb{P}\left( X_{n-k,n} \leq
		0\right),\mathbb{P}\left(\hat X_{n-k,n} \leq 0\right)\right) \to 0$ as
		well, and hence the Hill estimators based on either true observations or
		their approximations are well-defined with increasing probability. For
		the difference of the Hill estimators based on true and approximated
		observations (on the event where they are well-defined),
		\eqref{eq:ln_ineq} gives
		\begin{equation} \label{eq:gamma_diff_ineq}
			\begin{split}
				|\tilde\gamma_n - \hat\gamma_n|
				&\leq \frac{1}{k} \sum_{i=0}^{k-1}
				\left|\ln\left(\frac{\hat X_{n - i, n}}
				{X_{n - i, n}}\right)\right|
				+ \left|\ln\left(\frac{\hat X_{n - k, n}}{X_{n - k, n}}
				\right)\right|
				\leq \frac{2E_{n,n}}
				{\min\left(\hat X_{n-k, n}, X_{n-k, n}\right)} \\
				&= \frac{2E_{n,n}}{U(n/k)}
				\frac{U(n/k)}{\min\left(\hat X_{n-k, n}, X_{n-k, n}\right)}.
			\end{split}
		\end{equation}
		Here 
		 \begin{equation*}
			\begin{split}
				\left|\frac{\min\left(\hat X_{n-k, n}, X_{n-k, n}\right)}
				{U(n/k)} - 1\right|
				&\leq \frac{\left|\min\left(\hat X_{n-k, n}, X_{n-k, n}\right)
				- X_{n-k, n}\right|}{U(n/k)}
				+ \left|\frac{X_{n-k,n}}{U(n/k)} - 1\right| \\
				&\leq \frac{E_{n,n}}{U(n/k)}
				+ \left|\frac{X_{n-k,n}}{U(n/k)} - 1\right|,
			\end{split}
		\end{equation*}
		from which it follows that 
		\begin{equation} \label{eq:u_div_min}
			\frac{U(n/k)}{\min\left(\hat X_{n-k, n}, X_{n-k, n}\right)}
			\stackrel{\mathbb{P}}{\to} 1, \quad n\to\infty.
		\end{equation}
		For a fixed $\varepsilon>0$ set 
		\begin{equation*}
			A_n = \left\{|\tilde\gamma_n - \hat\gamma_n|
		  	\leq \frac{2E_{n,n}}{U(n/k)}\left(1 + \varepsilon\right)\right\}.
		\end{equation*}
		By \eqref{eq:gamma_diff_ineq} and~\eqref{eq:u_div_min} we have
		$\mathbb{P}\left(A_n^c\right) \to 0$. Since, for a positive sequence
		$c_n \in \mathbb{R}_{>0}$, 
		\begin{equation*}
		  	\begin{split}
		  		\mathbb{P}\left(\left|c_n\left(\tilde\gamma_n
		  		- \hat\gamma_n\right)\right| > \varepsilon\right)
		  		&\leq \mathbb{P}\left(c_n\left|\tilde\gamma_n
		  		- \hat\gamma_n\right| > \varepsilon, A_n\right)
		  		+ \mathbb{P}\left(A_n^c\right) \\
		  		&\leq \mathbb{P}\left(c_n\frac{E_{n,n}}{U(n/k)}
		  		> \frac{\varepsilon}{2(1 + \varepsilon)}\right)
				+ \mathbb{P}\left(A_n^c\right),
		  	\end{split}
		\end{equation*}
		it follows that $c_n\left(\tilde\gamma_n -
		\hat\gamma_n\right)\stackrel{\mathbb{P}}{\to} 0$ as $n\to\infty$
		whenever $c_n\frac{E_{n,n}}{U(n/k)} \stackrel{\mathbb{P}}{\to} 0$.
		Writing $c_n\left(\hat\gamma_n - \gamma\right) = c_n\left(\hat\gamma_n -
		\tilde\gamma_n\right) + c_n\left(\tilde\gamma_n - \gamma\right)$ and
		setting $c_n \equiv 1$, Part~\ref{item:hill-consistency} of Theorem
		\ref{theorem:approx} follows by~\cite[Theorem 3.2.2]{deHaan2007} and
		Slutsky's lemma. Similarly, setting $c_n = \sqrt{k}$ and using
		~\cite[Theorem 3.2.5]{deHaan2007} together with Slutsky's Lemma gives
		Part~\ref{item:hill-normality}.
		  
		It remains to prove Part~\ref{item:hill-quantile}. We write
		\begin{equation*}
			\begin{split}
				\frac{\sqrt{k}}{\ln d_n}\left(\frac{\hat x_q}{x_q} - 1\right)
				= \underbrace{\frac{\sqrt{k}}{\ln d_n}\left(\frac{\tilde x_q}
					{x_q} - 1\right)}_{\text{I}}
				+ \underbrace{\frac{\sqrt{k}}{\ln d_n}\left(\frac{\hat x_q}
					{x_q} - \frac{\tilde x_q}{x_q}\right)}_{\text{II}},
			\end{split}
		\end{equation*}
		where Term I converges in distribution to $N\left(\frac{\lambda}{1 -
		\rho},\gamma^2\right)$ by~\cite[Theorem 4.3.8]{deHaan2007}. For Term II
		we split
	  	\begin{equation*}
			\textnormal{II} = \underbrace{\frac{U(n/k)d_n^{\gamma}}{x_q}}
			_{\textnormal{III}}
			\underbrace{\frac{\sqrt{k}}{\ln d_n}\left(\frac{\hat X_{n-k,n}
			d_n^{\hat\gamma_n - \gamma}}{U(n/k)} - \frac{X_{n-k,n}
			d_n^{\tilde\gamma_n - \gamma}}{U(n/k)}\right)}_{\textnormal{IV}}.
	  	\end{equation*}
		It follows from the assumptions and \cite[Remark B.3.15]{deHaan2007}
		that
 		\begin{equation*}
			\textnormal{III} = \frac{U\left(\frac{n}{k}\right)d_n^\gamma}
			{U\left(\frac{n}{k}d_n\right)} \to 1, \quad n\to\infty
		\end{equation*}
		and hence it suffices to show that, as $n\to\infty$, Term IV converges
		to zero in probability. We write
		\begin{equation*}
			\begin{split}
				\textnormal{IV} &= \underbrace{\frac{\sqrt{k}}{\ln d_n}
				\left(d^{\hat\gamma_n - \gamma} - 1\right)
				\left(\frac{\hat X_{n-k,n}}{U(n/k)} - 1\right)}
				_{\textnormal{IV}_1}
				- \underbrace{\frac{\sqrt{k}}{\ln d_n}
				\left(d^{\tilde\gamma_n - \gamma} - 1\right)
				\left(\frac{X_{n-k,n}}{U(n/k)} - 1\right)}
				_{\textnormal{IV}_2} \\
				&+ \underbrace{\frac{\sqrt{k}}{\ln d_n}
				\left(d_n^{\hat\gamma_n - \gamma}
				- d_n^{\tilde\gamma_n - \gamma}\right)}_{\textnormal{IV}_3}
				+ \underbrace{\frac{\sqrt{k}}{\ln d_n}
				\frac{\hat X_{n-k,n} - X_{n-k,n}}{U(n/k)}}_{\textnormal{IV}_4}
			\end{split}
		\end{equation*}
		and show that each term converge. For this, note first that $1/\ln d_n
		\to 0$ since $d_n\to \infty$. As $\ln(nq) = o\left(\sqrt{k}\right)$ and
		$\sqrt{k}(\hat\gamma_n - \gamma) = O_\mathbb{P}(1)$, we also have
  		\begin{equation*}
			d_n^{\hat\gamma_n - \gamma} = \exp\left(\sqrt{k}
			(\hat\gamma_n - \gamma)
			\left(\frac{\ln(k)}{\sqrt{k}} - \frac{\ln(nq)}{\sqrt{k}}\right)
			\right) \stackrel{\mathbb{P}}{\to} 1, \quad n\to\infty.
		\end{equation*}
		Using the same reasoning gives $d_n^{\tilde\gamma_n - \gamma}
		\stackrel{\mathbb{P}}{\to} 1$.  

		Since $\sqrt{k} \left(\frac{X_{n-k, n}}{U(n/k)} - 1\right)$ is bounded
		in probability (see \cite[Theorem 4.3.8]{deHaan2007}) and since, by
		\eqref{eq:order_stat_diff_ineq} and assumptions, 
		\begin{equation} \label{eq:IV4}
			\sqrt{k} \frac{|\hat X_{n-k, n} - X_{n-k, n}|}{U(n/k)}
			\leq \sqrt{k}\frac{E_{n,n}}{U(n/k)}\stackrel{\mathbb{P}}{\to} 0,
		\end{equation}
		we have 
		\begin{equation} \label{eq:hat_bounded}
			\sqrt{k} \left(\frac{\hat X_{n-k, n}}{U(n/k)} - 1\right)
			= \sqrt{k} \left(\frac{X_{n-k, n}}{U(n/k)} - 1\right)
			+ \sqrt{k} \frac{\hat X_{n-k, n} - X_{n-k, n}}{U(n/k)}
			= O_{\mathbb{P}}(1).
		\end{equation}
		Now \eqref{eq:IV4} gives $\textnormal{IV}_4\stackrel{\mathbb{P}}{\to}0$
		and $\textnormal{IV}_1\stackrel{\mathbb{P}}{\to}0$ follows from
		\eqref{eq:hat_bounded}. Moreover, similar argumentation shows that
		$\textnormal{IV}_2\stackrel{\mathbb{P}}{\to}0$.

		For Term $\textnormal{IV}_3$, we have
  		\begin{equation*}
			\textnormal{IV}_3 = \frac{\sqrt{k}}{\ln d_n}
			\left(d_n^{\hat\gamma_n - \tilde\gamma_n} - 1\right)
			\left(d_n^{\tilde\gamma_n - \gamma} - 1\right)
			+ \frac{\sqrt{k}}{\ln d_n}
			\left(d_n^{\hat\gamma_n - \tilde\gamma_n} - 1\right)
		\end{equation*}
		and hence it suffices to show $\frac{\sqrt{k}}{\ln d_n}
		\left(d_n^{\hat\gamma_n - \tilde\gamma_n} - 1\right)
		\stackrel{\mathbb{P}}{\to} 0$. We have 
		\begin{equation*}
			\frac{\sqrt{k}}{\ln d_n}
			\left(d_n^{\hat\gamma_n - \tilde\gamma_n} - 1\right)
			= \sqrt{k}\left(\hat\gamma_n - \tilde\gamma_n\right)
			\frac{\exp\left(\left(\hat\gamma_n
			- \tilde\gamma_n\right)\ln d_n\right)- 1}
			{\left(\hat\gamma_n - \tilde\gamma_n\right)\ln d_n},
		\end{equation*}
		where $\sqrt{k}\left(\hat\gamma_n - \tilde\gamma_n\right)
		\stackrel{\mathbb{P}}{\to}0$ by the proof of
		Part~\ref{item:hill-normality}. This also implies 
		\begin{equation*}
			Y_n := \left(\hat\gamma_n - \tilde\gamma_n\right)\ln d_n
			= \sqrt{k}(\hat\gamma_n - \tilde\gamma_n)
			\left(\frac{\ln(k)}{\sqrt{k}} - \frac{\ln(nq)}{\sqrt{k}}\right)
			\stackrel{\mathbb{P}}{\to} 0.
		\end{equation*}
		Using the Taylor expansion for $e^x$ on the event $\{|Y_n|<1\}$ gives 
		\begin{equation*}
			\left|\frac{e^{Y_n} - 1}{Y_n} - 1\right|
			= \left|\sum_{j=1}^\infty \frac{Y_n^j}{(j+1)!}\right|
			\leq \sum_{j=1}^\infty \left|Y_n\right|^j
			= \frac{1}{1 - \left|Y_n\right|} - 1 \stackrel{\mathbb{P}}{\to} 0.
		\end{equation*}
		This implies $\textnormal{IV}_3\stackrel{\mathbb{P}}{\to}0$. 

		Combining all the above shows that, as $n\to \infty$, we have
		$\textnormal{IV}\stackrel{\mathbb{P}}{\to}0$. This completes the proof
		of Part~\ref{item:hill-quantile}.
	\end{proof}

	\begin{proof}[Proof of Theorem \ref{theorem:concentration}]
		Suppose that the assumptions of Theorem \ref{theorem:concentration}
		hold. We divide the proof into two steps. In the first step we consider
		the probability $\mathbb{P}\left(\left|\frac{X_{n-k, n}}{U(n/k)} -
		1\right| > x \right)$ and in the second step we consider the probability
		$\mathbb{P}\left(\left|\frac{U(n/k)}{X_{n-k, n}} - 1\right| > x
		\right)$. We denote $T = 1/(1 - F_X)$ so that  $U = T^\leftarrow$.
		Throughout the proof, let $Y_1, \ldots, Y_n$ be an i.i.d.\ sample from
		the standard exponential distribution $F_Y(x) = 1 - e^{-x}$, $x \geq 0$.
		Then $X_i \stackrel{d}{=} U\left(\exp\left(Y_i\right)\right)$.
		
		\paragraph*{Step 1.} Let first $0 < x < 1$ and $n$ large enough so that
		$U(n/k) > 0$. By using right-continuity of $T$ we obtain
		\begin{equation*}
			\begin{split}
				\mathbb{P}\left(\left|\frac{X_{n-k, n}}{U(n/k)} - 1\right| > x
				\right) &=
				\begin{multlined}[t]
					\mathbb{P}\left(U\left(\exp\left(Y_{n-k, n}\right)\right) >
					U\left(\frac{n}{k}\right)(1 + x)\right) \\
					+ \mathbb{P}\left(U\left(\exp\left(Y_{n-k, n}\right)\right)
					< U\left(\frac{n}{k}\right)(1 - x)\right)
				\end{multlined} \\
				&\leq
				\begin{multlined}[t]
					\mathbb{P}\left(\exp\left(Y_{n-k, n}\right) >
					T\left(U\left(\frac{n}{k}\right)(1 + x)\right)\right) \\
					+ \mathbb{P}\left(\exp\left(Y_{n-k, n}\right)
					\leq T\left(U\left(\frac{n}{k}\right)(1 - x)\right)\right).
				\end{multlined}
			\end{split}
		\end{equation*}
		Let $0 < \delta < \min\left\{(1+x)^{1/\gamma} - 1, 1 - (1-x)^{1/\gamma},
		\frac{1}{k} \left(1-x\right)^{1/\gamma}\right\}$. By Part
		\ref{item:F-RV} of Theorem \ref{theorem:heavy} we have, for all $y > 0$,
		that $T(U(n/k)y)/T(U(n/k)) \to y^{1/\gamma}$ as $n\to\infty$. Moreover,
		$T\left(U\left(n/k\right)\right) \sim n/k$, see \cite[Section
		2.2.1]{resnick2007}. Thus, there exists $N=N\left(\delta, y\right)$ such
		that for $n\geq N$ we have 
		\begin{equation*}
			\frac{n}{k}\left(y^{1/\gamma} - \delta\right) 
			\leq T\left(U(n/k)y\right) 
			\leq \frac{n}{k}\left(y^{1/\gamma} + \delta\right),
			\quad y > 0.
		\end{equation*}
		Set $M_1 = (1+x)^{1/\gamma} - \delta$ and $M_2 = (1-x)^{1/\gamma} +
		\delta$. Then, for large enough $n$, 
		\begin{equation*}
			\mathbb{P}\left(\left|\frac{X_{n-k, n}}{U(n/k)} - 1\right| > x
			\right) \leq \underbrace{\mathbb{P}\left(\exp\left(Y_{n-k, n}\right)
			> M_1 \frac{n}{k}\right)}_{\textnormal{I}}
			+ \underbrace{\mathbb{P}\left(\exp\left(Y_{n-k, n}\right)
			\leq M_2 \frac{n}{k}\right)}_{\textnormal{II}}.
		\end{equation*}
		We next bound the terms by using the moment generating function of
		$Y_{n-k,n}$, denoted by $g(a)$. By R\'enyi's representation we have
		\begin{equation*}
			Y_{n-k, n} \stackrel{d}{=} \sum_{j=1}^{n-k} \frac{E_j}{n-j+1},
		\end{equation*}
		where $E_1, \ldots, E_{n-k}$ are i.i.d.\ random variables from the
	 	exponential distribution $\text{Exp}(1)$. This gives
		\begin{equation*}
			g(a) = \prod_{j=1}^{n-k} g_{E_1}\left(\frac{a}{n-j+1}\right)
			= \prod_{j=1}^{n-k} \left(1 - \frac{a}{n-j+1}\right)^{-1}
			= \prod_{j=k+1}^{n} \frac{j}{j-a},
		\end{equation*}
		provided that $a < k+1$, where $g_{E_1}$ is the moment generating
		function of $\text{Exp}(1)$. 
		
		Consider now Term I and let $0 < a < k + 1$. By the Markov inequality
		\begin{equation*}
			\begin{split}
				\textnormal{I} &= \mathbb{P}\left(\exp\left(aY_{n-k, n}\right) >
				\left(M_1 \frac{n}{k}\right)^a \right)
				\leq g(a) \left(\frac{nM_1}{k}\right)^{-a}
				= \left(\frac{k}{nM_1}\right)^{a}
				\prod_{j=k+1}^{n} \frac{j}{j-a} \\
				&= \left(\frac{k}{nM_1}\right)^{a}
				\frac{\Gamma\left(n+1\right)\Gamma\left(k+1-a\right)}
				{\Gamma\left(k+1\right)\Gamma\left(n+1-a\right)}
				=: f(a),
			\end{split}
		\end{equation*}
		where $\Gamma$ is the Gamma function. Our aim is to choose $a$ such that
		$f(a)$, or equivalently $\ln f(a)$, is minimized. Differentiating gives 
		\begin{equation*}
			\begin{split}
				\frac{\mathrm{d} \ln f(a)}{\mathrm{d}a}
				&= \ln\left(\frac{k}{nM_1}\right)
				+ \sum_{j=k+1}^n \frac{1}{j-a} \\
				&= \left(\psi\left(n+1-a\right) - \ln n\right)
				- \left(\psi\left(k+1-a\right) - \ln k + \ln M_1\right),
			\end{split}
		\end{equation*}
		where $\psi$ is the Digamma function, and 
			\begin{equation*}
			\frac{\mathrm{d}^2 \ln f(a)}{\mathrm{d}a^2}
			= \sum_{j=k+1}^{n} \frac{1}{(j-a)^2} > 0.
		\end{equation*}
		Thus $\ln f(a)$ is convex, and hence the unique minimum is obtained at
		$a_1$ solving $\frac{\mathrm{d} \ln f(a)}{\mathrm{d}a} = 0$. We next
		examine the asymptotic behavior of $a_1$. For this we use the asymptotic
		formulae
		\begin{align}
			\psi(z) &= \ln z + O\left(\frac{1}{z}\right), \quad z\to\infty,
			\label{eq:digamma} \\
			\ln\left(1 + z\right) &= z + o\left(z\right), \quad z\to 0
			\label{eq:log}
		\end{align}
		leading, as $n\to\infty$, to
		\begin{align*}
			\frac{\mathrm{d} \ln f(k)}{\mathrm{d}a} &= \ln k + O(1) > 0
			\quad\textnormal{and} \\
			\frac{\mathrm{d} \ln f(1)}{\mathrm{d}a} &= -\ln M_1 + o(1) < 0,
		\end{align*}
		since $M_1 > 1$. It follows that $1\leq a_1 \leq k$ and hence $n + 1 -
		a_1\to\infty$. Denote $b_n' = \psi\left(n+1-a_1\right) - \ln n$. Now, by
		\eqref{eq:digamma} and~\eqref{eq:log} we have
		\begin{equation*}
			b'_n = \ln\left(1 + \frac{1 - a_1}{n}\right) + o(1)
			= \frac{1-a_1}{n} + o(1) = o(1).
		\end{equation*}
		Since $a_1$ solves $\frac{\mathrm{d} \ln f(a)}{\mathrm{d}a} = 0$, it
		follows that
		\begin{equation*}
			\psi\left(k+1-a_1\right) - \ln k + \ln M_1 = o(1).
		\end{equation*}
		Since now $M_1$ is fixed, this further implies that $\psi\left(k + 1 -
		a_1\right) - \ln k = O(1)$ and hence $\psi\left(k + 1 -
		a_1\right)\to\infty$. From monotonicity of $\psi$ we also get $k + 1 -
		a_1\to\infty$ allowing to apply \eqref{eq:digamma} that leads to
		\begin{equation*}
			\ln\left(\frac{M_1(k + 1 - a_1)}{k}\right) = o(1),
		\end{equation*}
		and thus, as $n\to\infty$, the minimizer $a_1$ satisfies
		\begin{equation*}
			a_1 = k + 1 - \frac{k}{M_1} + b_n,
		\end{equation*}
		where $b_n = o\left(k\right)$. We now evaluate the minimum value
		$f(a_1)$. Stirling's approximation~$ \Gamma(x + 1) \sim \sqrt{2\pi x}
		\left(\frac{x}{e}\right)^x$ gives, for any sequence $a_n\in\mathbb{R}$
		such that $n+1+a_n\to\infty$ and $k+1+a_n\to\infty$, that
		\begin{equation*}
			\frac{\Gamma(n + 1)\Gamma(k + 1 + a_n)}
			{\Gamma(n + 1 + a_n)\Gamma(k + 1)}
			\sim \left(1 + \frac{a_n}{k} \right)^{k + a_n + \frac{1}{2}}
			\left(1 + \frac{a_n}{n}\right)^{-(n + a_n + \frac{1}{2})}
			k^{a_n}n^{-a_n}, \quad n\to\infty.
		\end{equation*}
		This leads to 
		\begin{equation*}
			f(a_1) \sim \left(\frac{1}{M_1}\right)^{a_1}
			\left(1 - \frac{a_1}{k} \right)^{k - a_1 + \frac{1}{2}}
			\left(1 - \frac{a_1}{n}\right)^{-(n - a_1 + \frac{1}{2})}.
		\end{equation*}
		On the other hand, \eqref{eq:log} implies that
		\begin{equation*}
			\left(1 + \frac{y}{z}\right)^{z(1+o(1))}
			= \exp\left(y + o(y)\right)
			\quad\textnormal{if}\quad y = o(z), \quad z\to\infty.
		\end{equation*}
		This gives 
		\begin{equation*}
			\begin{split}
				\left(1 - \frac{a_1}{k} \right)^{k - a_1 + \frac{1}{2}} =& 
				\left(\frac{1}{M_1}\right)^{\frac{k}{M_1}-\frac{1}{2}-b_n}
				\left(1 - \frac{1 + b_n}{k/M_1}\right)^{\frac{k}{M_1}(1 + o(1))}
				\\
				=& 
				\left(\frac{1}{M_1}\right)^{\frac{k}{M_1}-\frac{1}{2}-b_n}
				\exp\left(-b_n - 1 + o(k)\right)
			\end{split}
		\end{equation*} 
		and
		\begin{equation*}
			\left(1 - \frac{a_1}{n}\right)^{-(n - a_1 + \frac{1}{2})}
			= \left(1 - \frac{a_1}{n}\right)^{-n(1 + o(1))}
			= \exp(a_1 + o(k)), \quad n\to\infty.
		\end{equation*}
		Hence, we acquire
		\begin{equation*}
			f(a_1) \sim \exp\left(k - \frac{k}{M_1} + o(k)\right)
			\left(\frac{1}{M_1}\right)^{k+\frac{1}{2}}
		\end{equation*}
		that can be equivalently written as
		\begin{equation*}
			f(a_1) = (1 + o(1))\left(\frac{1}{M_1}\right)^{1/2}
			\exp\left(k - \frac{k}{M_1} - k\ln M_1 + o(k)\right),
			\quad n\to\infty.
		\end{equation*}
		Here, since $\delta < (1+x)^{1/\gamma} - 1$, we have $M_1^{-1/2} < 1$.
		Recall also that $\textnormal{I} \leq f(a)$. Thus we can find $N =
		N(\delta,x)$ such that, for $n\geq N$, we have
		\begin{equation*}
			\textnormal{I} \leq \left(1 + \delta\right)
			\exp\left(\left(1 - \frac{1}{M_1} - \ln M_1 + \delta\right)k\right).
		\end{equation*}
		It remains to write this upper bound in terms of $x$. For this we use
		Maclaurin series of $\ln(1-z)$ to get
		\begin{equation*}
			\ln\left(1 - z\right) \geq -\frac{z}{1-z}, \quad 0 < z < 1.
		\end{equation*}
		Thus, since $\delta < \left(1 + x\right)^{1/\gamma} - 1$,
		\begin{equation*}
			\begin{split}
				\ln M_1 &= \frac{\ln\left(1 + x\right)}{\gamma}
				+ \ln\left(1 - \frac{\delta}{(1+x)^{1/\gamma}}\right)
				\geq \frac{\ln\left(1 + x\right)}{\gamma}
				- \frac{\delta}{(1+x)^{1/\gamma} - \delta} \\
				&\geq \frac{\ln\left(1 + x\right)}{\gamma} - \delta.
			\end{split}
		\end{equation*}
		Together with $\frac{1}{M_1} \geq (1+x)^{-1/\gamma}$ this leads to 
		\begin{equation} \label{eq:I}
			\textnormal{I} \leq \left(1 + \delta\right)
			\exp\left(\left(1 - (1 + x)^{-1/\gamma}
			- \frac{\ln\left(1 + x\right)}{\gamma} + 2\delta\right)k\right).
		\end{equation}
		
		We next treat Term II with similar arguments. Let $a > 0$. By the Markov
		inequality
		\begin{equation*}
			\begin{split}
				\textnormal{II} &= \mathbb{P}\left(\exp\left(-aY_{n-k, n}\right)
				\geq \left(M_2 \frac{n}{k}\right)^{-a} \right)
				\leq g(-a) \left(\frac{nM_2}{k}\right)^{a}
				= \left(\frac{nM_2}{k}\right)^{a}
				\prod_{j=k+1}^{n} \frac{j}{j+a} \\
				&= \left(\frac{nM_2}{k}\right)^{a}
				\frac{\Gamma\left(n+1\right)\Gamma\left(k+1+a\right)}
				{\Gamma\left(k+1\right)\Gamma\left(n+1+a\right)}
				=: h(a).
			\end{split}
		\end{equation*}
		As in bounding Term I, we choose $a$ that minimizes $h(a)$. We obtain
		\begin{equation*}
			\begin{split}
				\frac{\mathrm{d} \ln h(a)}{\mathrm{d}a}
				&= \ln\left(\frac{nM_2}{k}\right)
				- \sum_{j=k+1}^n \frac{1}{j+a} \\
				&= \left(\psi\left(k+1+a\right) - \ln k + \ln M_2\right)
				- \left(\psi\left(n+1+a\right) - \ln n \right)
			\end{split}
		\end{equation*}
		and
		\begin{equation*}
			\frac{\mathrm{d}^2 \ln h(a)}{\mathrm{d}a^2}
			= \sum_{j=k+1}^{n} \frac{1}{(j+a)^2} > 0
		\end{equation*}
		giving us unique solution $a_2$ to $\frac{\mathrm{d} \ln
		h(a)}{\mathrm{d}a} = 0$. Recall that $M_2 = (1-x)^{1/\gamma} + \delta$,
		where $\delta < 1 - (1-x)^{1/\gamma}$. Now, as $n\to \infty$,
		\eqref{eq:digamma} and~\eqref{eq:log} give
		\begin{align*}
			\frac{\mathrm{d} \ln h(1)}{\mathrm{d}a} &= \ln M_2 < 0
			\quad\textnormal{and} \\
			\frac{\mathrm{d} \ln h\left(\frac{k}{M_2}\right)}{\mathrm{d}a}
			&= \ln\left(1 + M_2\right) + o(1) > 0.
		\end{align*}
		It follows that $1\leq a_2 \leq \frac{k}{M_2}$. Completing similar steps
		as in the case of Term I gives us, as $n\to \infty$, that
		\begin{equation*}
			a_2 = \frac{k}{M_2} - k - 1 + o(k)
		\end{equation*}
		and
		\begin{equation*}
			h(a_2) = (1 + o(1)) \left(\frac{1}{M_2}\right)^{1/2}
			\exp\left(k - \frac{k}{M_2} - k\ln M_2 + o(k)\right).
		\end{equation*}
		Note that in this case $M_2^{-\frac{1}{2}} > 0$, and hence we obtain the
		bound
		\begin{equation*}
			\textnormal{II} \leq h(a_2) = (1 + o(1)) 
			\exp\left(k - \frac{k}{M_2} - (k+1/2)\ln M_2 + o(k)\right).
		\end{equation*}
		In order to write this in terms of $x$, we choose $n$ large enough such
		that
		\begin{equation} \label{eq:x_bound}
			\frac{k+1}{k}(1-x)^{1/\gamma} < \frac{k}{k+1/2}.
		\end{equation}
		Note that
		\begin{equation*}
			\tilde b(a) := \exp\left(k - \frac{k}{a} - (k+1/2)\ln a
			\right)
		\end{equation*}
		is increasing for $a < \frac{k}{k+1/2}$. Together
		with~\eqref{eq:x_bound}, as $\delta < \frac{1}{k}(1-x)^{1/\gamma}$,
		this gives
		\begin{equation*}
			\textnormal{II} \leq (1 + o(1)) 
			\exp\left(k - \frac{k (1-x)^{-1/\gamma}}{1+\frac{1}{k}}
			- \left(k+\frac{1}{2} \right)\ln \left( \left(1+\frac{1}{k}\right)
			(1-x)^{1/\gamma} \right) + o(k)\right).
		\end{equation*}
		That is, we can find $N=N(\delta,x)$ such that, for $n\geq N$, we have
		\begin{equation}
			\label{eq:II}
			\textnormal{II} \leq 
			(1 + \delta)(1-x)^{-1/(2\gamma)} b_{\gamma, \delta, k}(1-x)
			\exp\left((1-x)^{-1/\gamma} \frac{k}{k+1}\right).
		\end{equation}
		By choosing $\varepsilon = 2\delta$ and combining~\eqref{eq:I} and
		\eqref{eq:II} leads to
		\begin{align*}
			\mathbb{P}\left(\left|\frac{X_{n-k, n}}{U(n/k)} - 1\right| > x
			\right) \leq & \left(1 + \varepsilon\right)\Bigg(
			b_{\gamma, \varepsilon, k}\left(1+x\right) \\
			& + (1-x)^{-\frac{1}{2\gamma}} b_{\gamma, \varepsilon, k}
			\left(1 - x\right)\exp \left((1-x)^{-1/\gamma}
			\frac{k}{k+1}\right)\Bigg).
		\end{align*}	
		This proves the claimed concentration for $0 < x < 1$. 

		Suppose next $x \geq 1$. Then 
		\begin{equation*}
			\begin{split}
				\mathbb{P}\left(\left|\frac{X_{n-k, n}}{U(n/k)} - 1\right| > x
				\right)
				&= \mathbb{P}\left(X_{n-k,n} > U\left(\frac{n}{k}\right)
				\left(1 + x\right)\right)
				+ \mathbb{P}\left(X_{n-k,n} < U\left(\frac{n}{k}\right)
				\left(1 - x\right)\right) \\
				&\leq \mathbb{P}\left(X_{n-k,n} > U\left(\frac{n}{k}
				\right)\left(1 + x\right)\right)
				+ \mathbb{P}\left(X_{n-k,n} < 0\right).
			\end{split}
		\end{equation*}
		Here the first term can be treated as Term I in the case $0 < x < 1$,
		leading to
		\begin{equation*}
			\mathbb{P}\left(\left|\frac{X_{n-k, n}}{U(n/k)} - 1\right| > x
			\right)
			\leq \left(1 + \varepsilon\right)
			b_{\gamma, \varepsilon, k}\left(1+x\right)
			+ \mathbb{P}\left(X_{n-k,n} < 0\right)
		\end{equation*}
		for $n\geq N=N(\varepsilon,x)$. 
		This completes Step 1 of the proof.

		\paragraph*{Step 2.} Suppose first that $0 < x < 1$. Then
		\begin{align*}
			\mathbb{P}\left(\left|\frac{U(n/k)}{X_{n-k, n}} - 1\right| > x
			\right)
			\leq& \mathbb{P}\left(\frac{U(n/k)}{X_{n-k, n}} > 1 + x,
			X_{n-k,n} > 0\right)
			+ \mathbb{P}\left(\frac{U(n/k)}{X_{n-k, n}} < 1 - x,
			X_{n-k,n} > 0\right) \\
			&+ \mathbb{P}\left(X_{n-k,n} \leq 0\right).
		\end{align*}
		For the first term we have
		\begin{equation*}
			\mathbb{P}\left(\frac{U(n/k)}{X_{n-k, n}} > 1 + x, X_{n-k,n} > 0
			\right)
			\leq \mathbb{P}\left(X_{n-k,n}
			< U(n/k)\left(1 - \frac{x}{1+x}\right)\right)
		\end{equation*}
		and for the second term we have
		\begin{equation*}
			\mathbb{P}\left(\frac{U(n/k)}{X_{n-k, n}} < 1 - x, X_{n-k,n} > 0
			\right)
			\leq \mathbb{P}\left(X_{n-k,n} > U(n/k)\left(1 + \frac{x}{1-x}
			\right)\right),
		\end{equation*}
		where $\mathbb{P}\left(X_{n-k,n} < U(n/k)\left(1 -
		\frac{x}{1+x}\right)\right)$ and $\mathbb{P}\left(X_{n-k,n} >
		U(n/k)\left(1 + \frac{x}{1-x}\right)\right)$ can be handled similarly as
		in Step 1. This gives, for $0 < x < 1$, that 
		\begin{equation*}
			\begin{split}
				\mathbb{P}\left(\left|\frac{U(n/k)}{X_{n-k, n}} - 1\right| > x
				\right)
				\leq& \left(1 + \varepsilon\right)(1+x)^{\frac{1}{2\gamma}}
				b_{\gamma, \varepsilon, k}\left(\frac{1}{1+x}\right)
				\exp \left((1+x)^{1/\gamma} \frac{k}{k+1}\right) \\
				&+ \left(1 + \varepsilon\right)
				b_{\gamma, \varepsilon, k}\left(\frac{1}{1-x}\right)
				+ \mathbb{P}\left(X_{n-k,n} \leq 0\right)
			\end{split}
		\end{equation*}
		for $n\geq N = N(\varepsilon,x)$. Similarly for $x\geq 1$,
		\begin{equation*}
			\begin{split}
				\mathbb{P}\left(\left|\frac{U(n/k)}{X_{n-k, n}} - 1\right| > x
				\right)
				\leq \mathbb{P}\left(\frac{U(n/k)}{X_{n-k, n}} > 1 + x,
				X_{n-k,n} > 0\right)
				+ \mathbb{P}\left(X_{n-k,n} \leq 0\right) \\
				\leq \mathbb{P}\left(X_{n-k,n}
				< U(n/k)\left(1 - \frac{x}{1+x}\right)\right)
				+ \mathbb{P}\left(X_{n-k,n} \leq 0\right),
			\end{split}
		\end{equation*}
		where $\mathbb{P}\left(X_{n-k,n} < U(n/k)\left(1 -
		\frac{x}{1+x}\right)\right)$ can be handled as in Step 1. This gives,
		for $x \geq 1$, that 
		\begin{align*}
			\mathbb{P}\left(\left|\frac{U(n/k)}{X_{n-k, n}} - 1\right| > x
			\right)
			\leq& \left(1 + \varepsilon\right)(1+x)^{\frac{1}{2\gamma}}
			b_{\gamma, \varepsilon, k}\left(\frac{1}{1+x}\right)
			\exp \left((1+x)^{1/\gamma}\frac{k}{k+1}\right) \\
			&+ \mathbb{P}\left(X_{n-k,n} \leq 0\right)
		\end{align*}
		for $n\geq N=N(\varepsilon,x)$. This completes Step 2.
	\end{proof}

	The following proposition, required for the proofs of Theorems
	\ref{theorem:continuity} and \ref{theorem:monte}, follows from Theorem
	\ref{theorem:concentration}.
	\begin{proposition} \label{prop:min-concentration} Let $X_1, \ldots, X_n$ be
		i.i.d.\ copies of almost surely positive $X\in RV_\gamma$, and $\hat
		X_1, \ldots, \hat X_n$ the corresponding almost surely positive
		approximations. Assume that, as $n\to \infty$, $k = k_n\to\infty$ and
		$k/n\to 0$. Let $\varepsilon>0$ be fixed and let $b_{\gamma,
		\varepsilon, k}\left(x\right)$ be given by \eqref{eq:b}. 
		\begin{enumerate}
			\item Suppose $0 < y < 1$. Then for any $0 < z < \frac{y}{1+y}$
			there exists $N =N(\varepsilon,y,z)$ such that for $n\geq N$ we
			have
			\footnotesize
			\begin{multline*}
				\mathbb{P}\left(\left|\frac{U(n/k)}
				{\min\left(\hat X_{n-k,n}, X_{n-k,n}\right)}-1\right| > y\right)
				\leq \left(1 + \varepsilon\right)
				\bigg(b_{\gamma, \varepsilon, k}\left(\frac{1}{1-y}
				- z\right) 
				\\
				+  \left(\frac{1}{1+y}+z\right)^{-1/(2\gamma)}
				b_{\gamma, \varepsilon, k}\left(\frac{1}{1+y} + z\right)
				\exp \left(\left(\frac{1}{1+y} + z \right)^{-1/\gamma}
				\frac{k}{k+1}  \right)
				\bigg)
				\\
				+ \mathbb{P}\left(\frac{E_{n,n}}{U(n/k)} > z\right).
			\end{multline*}
			\item Suppose $y \geq 1$. Then for any $0 < z < \frac{y}{1+y}$
			there exists $N =N(\varepsilon,y,z)$ such that for $n\geq N$ we
			have
			\footnotesize
			\begin{multline*}
				\mathbb{P}\left(\left|\frac{U(n/k)}
				{\min\left(\hat X_{n-k,n}, X_{n-k,n}\right)}-1\right| > y\right)
				\\
				\leq \left(1 + \varepsilon\right)
				\left(\frac{1}{1+y}+z\right)^{-1/(2\gamma)}
				b_{\gamma, \varepsilon, k}\left(\frac{1}{1+y} + z\right)
				\exp\left( \left(\frac{1}{1+y} + z \right)^{-1/\gamma}
				\frac{k}{k+1}  \right) \\
				+ \mathbb{P}\left(\frac{E_{n,n}}{U(n/k)} > z\right).
			\end{multline*}
		\end{enumerate}
	\end{proposition}

	\begin{proof}
		Denote $\tilde X_n = \min\left(\hat X_{n-k, n}, X_{n-k, n}\right)$ and
		$\tilde E_n = \tilde X_n - X_{n-k,n}$. Since $|\hat X_{n-k,n} -
		X_{n-k,n}| \leq E_{n,n}$, it follows that $|\tilde E_n| \leq E_{n,n}$.
		Then
		\begin{equation*}
			\scriptsize
			\begin{split}
				\mathbb{P}\left(\left|\frac{U(n/k)}{\tilde X_n} - 1\right| > y
				\right)
				\leq& \mathbb{P}\left(\left|\frac{U(n/k)}{\tilde X_n} - 1\right|
				> y, \left|\tilde E_n\right| \leq z U(n/k)\right)
				+ \mathbb{P}\left(\left|\tilde E_n\right| > z U(n/k)\right) \\
				\leq& \underbrace{\mathbb{P}\left(\frac{U(n/k)}{\tilde X_n}
				< 1 - y,\left|\tilde E_n\right| \leq z U(n/k)\right)}
				_{\textnormal{I}}
				+ \underbrace{\mathbb{P}\left(\frac{U(n/k)}{\tilde X_n} > 1 + y,
				\left|\tilde E_n\right| \leq z U(n/k)\right)}
				_{\textnormal{II}} \\
				&+ \mathbb{P}\left(\frac{E_{n,n}}{U(n/k)} > z\right).
			\end{split}
		\end{equation*}
		For Term I in the case $y\geq 1$ and $U(n/k) > 0$, we get
		$\textnormal{I} \leq \mathbb{P}\left(\frac{U(n/k)}{\tilde X_n} <
		0\right) = \mathbb{P}\left(\tilde X_n < 0\right) = 0$. In the case $0 <
		y < 1$ we have 
		\begin{equation*}
			\begin{split}
				\textnormal{I} &= \mathbb{P}\left(\tilde X_n > U(n/k)\frac{1}
				{1-y}, \left|\tilde E_n\right| \leq z U(n/k)\right) \\
				&= \mathbb{P}\left(X_{n-k,n} + \tilde E_n > U(n/k)\frac{1}
				{1-y}, \left|\tilde E_n\right|
				\leq z U(n/k)\right) \\
				&\leq \mathbb{P}\left(X_{n-k,n}
				> U(n/k)\left(\frac{1}{1-y} - z\right)\right)\\
				&=\mathbb{P}\left(X_{n-k,n} > U(n/k)
				\left(1 + x^{(1)}\right)\right),
			\end{split}
		\end{equation*}
		where $x^{(1)} = \frac{y}{1-y} - z$. This corresponds to Term I in
		Step 1 of the proof of Theorem \ref{theorem:concentration}, and thus 
		\begin{equation*} 
			\textnormal{I} \leq
			\left(1 + \varepsilon\right)
			b_{\gamma, \varepsilon, k}\left(1 + x^{(1)}\right).
		\end{equation*}
		For Term II we have, for any $y>0$, that
		\begin{equation*}
			\begin{split}
				\textnormal{II}
				&= \mathbb{P}\left(X_{n-k,n} + \tilde E_n
				< U(n/k)\frac{1}{1 + y},
				\left|\tilde E_n\right| \leq z U(n/k)\right) \\
				&\leq \mathbb{P}\left(X_{n-k,n}
				< U(n/k)\left(\frac{1}{1 + y} + z\right)\right)\\
				&=\mathbb{P}\left(X_{n-k,n} < U(n/k)
				\left(1 - x^{(2)}\right)\right)
			\end{split}
		\end{equation*}
		with $x^{(2)} = \frac{y}{1+y} - z$. This corresponds to Term II in
		Step 1 of the proof of Theorem \ref{theorem:concentration}, and hence
		\begin{equation*} 
			\textnormal{II} \leq
			\left(1 + \varepsilon\right)\left(1 - x^{(2)}\right)
			^{-\frac{1}{2\gamma}}
			b_{\gamma, \varepsilon, k}\left(1 - x^{(2)}\right)
			\exp\left( \left(1-x^{(2)} \right)^{-1/\gamma} \frac{k}{k+1}\right).
		\end{equation*}
		This gives the result.
	\end{proof}

	\section{Proofs of Theorem \ref{theorem:continuity} and Theorem
	\ref{theorem:monte}}
	\label{appendix:sec:fda}
  
	The following lemma follows straightforwardly from \cite[Theorem
	5.3.12]{deHaan2007} for $\gamma > 0$. 
	\begin{lemma} \label{lemma:large_deviations} Let $X\in 2ERV_{\gamma,
		\tilde\rho}$ with $\gamma > 0$ and $\tilde\rho < 0$. Let $F$ and $F_G$
		denote the cumulative distribution functions of $X$ and the Fr\'echet
		distribution, respectively. Then, as $x_n\to \infty$, we have
		\begin{equation*}
			\lim_{n\to\infty} \frac{1 - F^n(U(n)x_n)}{1 - F_G(x_n)} = 1.
		\end{equation*}
	\end{lemma}

	\begin{proof}[Proof of Theorem \ref{theorem:continuity}]
	  	Suppose that the assumptions of Theorem \ref{theorem:continuity} hold.
	  	As a first step we show that
		\begin{equation} \label{eq:error_continuity}
	  		E_{n,n}
	  		= O_{\mathbb{P}}\left(\phi\left(\frac{1}{m}\right) U_V(n)\right).
	  	\end{equation}
	 	This, together with Theorem \ref{theorem:approx}, implies Parts
	 	\ref{part1-continuity}-\ref{part3-continuity}.
 
	  	Let
	  	\begin{equation*}
	  		\hat Y_i = \sum_{j = 0}^{m - 1} \mathbbm{1}_{A_j}
	  		Y_i\left(\frac{j}{m}\right),
	  	\end{equation*}
	  	where
	  	\begin{equation*}
	  		A_j =
	  		\begin{cases}
				[j/m, (j + 1) / m), & j \in\{0, \ldots, m - 2\} \\
				[j/m, (j + 1) / m], & j = m - 1,
	  		\end{cases}
	  	\end{equation*}
	  	so that $\hat X_i = \| \hat Y_i \|_p$. Now, for $1\leq p<\infty$,
	  	The Minkowski inequality and \eqref{eq:criterion} lead to
	  	\begin{equation*}
	  		\begin{split}
				|\left\|Y_i\right\|_p - \hat X_i|
				& \leq \left\|Y_i - \hat Y_i\right\|_p
				= \left(\sum_{j = 0}^{m - 1} \int_{A_j}
				\left|Y_i(t) - Y_i(j/m)\right|^p \,\mathrm{d}t\right)^{1/p} \\
				&\leq V_i \left(\sum_{j = 0}^{m - 1} \int_{A_j}
				\phi^p\left(\left|t - \frac{j}{m}\right|\right)
				\,\mathrm{d}t\right)^{1/p} \leq V_i \phi\left(\frac{1}{m}
				\right).
	  		\end{split}
	  	\end{equation*}
	  	Similarly, for $p=\infty$, we have $\left|\left\|Y_i\right\|_\infty -
	  	\hat X_i\right| \leq \left\|Y_i - \hat Y_i\right\|_\infty \leq  V_i
	  	\phi\left(\frac{1}{m}\right)$ and hence, for all $p\in [1,\infty]$, 
	  	\begin{equation} \label{eq:continuity}
	  		E_{n,n} = \max_{i\in\{1, \ldots, n\}}
			|\left\|Y_i\right\|_p - \hat X_i|
	  		\leq \phi\left(\frac{1}{m}\right) V_{n,n}.
	  	\end{equation}
	  	As $V_i\in RV_{\gamma'}$, \eqref{eq:error_continuity} now follows from
	  	\eqref{eq:continuity} and Theorem \ref{theorem:heavy}.
	  
	   	It remains to prove Part \ref{part4-continuity}. Let $x > 0$, $y \geq 1$
		and $0 < z < \frac{y}{1+y}$. Denote $X_i = \left\|Y_i\right\|_p$ and let
	  	\begin{equation*}
	  		A_n = \left\{\left|\frac{U_{\left\|Y\right\|_p}(n/k)}
	  		{\min\left(\hat X_{n-k, n}, X_{n-k, n}\right)} - 1\right|
	  		\leq y\right\}.
	  	\end{equation*}
	  	Using \eqref{eq:gamma_diff_ineq} now gives
	  	\begin{equation*}
	  		\mathbb{P}\left(\left|\hat\gamma_n  - \tilde\gamma_n\right| >
			x\right)
	  		\leq \underbrace{\mathbb{P}\left(A_n^c\right)}_{\textnormal{I}}
	  		+ \underbrace{\mathbb{P}\left(\frac{E_{n,n}}
	  		{\min\left(\hat X_{n-k, n}, X_{n-k, n}\right)} > \frac{x}{2},
	  		A_n\right)}_{\textnormal{II}}.
	  	\end{equation*}
		For Term I, Proposition \ref{prop:min-concentration} allows to estimate
 		\begin{equation*} 
 			\begin{split}
 				\textnormal{I} \leq \left(1 + \varepsilon\right)
 				&\left(a_{y,z}\right)^{-1/(2\gamma)}
 				b_{\gamma, \varepsilon, k}\left(a_{y,z}\right)
				\exp\left(\left(a_{y,z} \right)^{-1/\gamma}
				\frac{k}{k+1}\right) \\
 				&+\mathbb{P}\left(\frac{E_{n,n}}{U_{\|Y\|_p}(n/k)} > z
				\right).
 			\end{split}
 		\end{equation*}
		Since $e_n \to 0$ as $n\to\infty$ and since $1 - F_G(x) \sim
		\left(\gamma'x\right)^{-\frac{1}{\gamma'}}$, Lemma
		\ref{lemma:large_deviations} and \eqref{eq:continuity} give that, for
		$n\geq N$,
 		\begin{equation*}
 			\begin{split}
 				\mathbb{P}\left(\frac{E_{n,n}}{U_{\|Y\|_p}(n/k)} > z\right)
 				&\leq \mathbb{P}\left(V_{n,n} > U_V(n) z
 				\frac{U_{\left\|Y\right\|_p}(n/k)}
 				{\phi\left(\frac{1}{m}\right)U_V(n)}\right)
 				= 1 - F_V^n\left(U_V(n) \frac{z}{e_n}\right) \\
 				&\sim 1 - F_G\left(\frac{z}{e_n}\right)
				\sim \left(\frac{e_n}{\gamma'z}\right)^{1/\gamma'}
 			\end{split}
 		\end{equation*}
		providing
 		\begin{equation*}
			\mathbb{P}\left(\frac{E_{n,n}}{U_{\|Y\|_p}(n/k)} > z\right)
 			\leq \left(1 + \varepsilon\right)
 			\left(\frac{e_n}{\gamma'z}\right)^{1/\gamma'}.
 		\end{equation*}
 		For Term II, we have
		\begin{equation*}
			\begin{split}
				\textnormal{II}
				&= \mathbb{P}\left(\frac{E_{n,n}}{U_{\left\|Y\right\|_p}(n/k)}
				\frac{U_{\left\|Y\right\|_p}(n/k)}
				{\min\left(\hat X_{n-k, n}, X_{n-k, n}\right)} > \frac{x}{2},
				1 - y \leq \frac{U_{\left\|Y\right\|_p}(n/k)}
				{\min\left(\hat X_{n-k, n}, X_{n-k, n}\right)}
				\leq 1 + y\right) \\
				&\leq \mathbb{P}\left(\frac{E_{n,n}}
				{U_{\left\|Y\right\|_p}(n/k)} > \frac{x}{2\left(1+y\right)}
				\right).
			\end{split} 
		\end{equation*}
		Now, for $n\geq N$,
 		\begin{equation*} 
			\textnormal{II} \leq \left(1 + \varepsilon\right)
			\left(\frac{2\left(1 + y\right)e_n}{\gamma'x}\right)^{1/\gamma'}.
		\end{equation*}
		Combining the above estimates completes the proof.
	\end{proof}

	The following lemma follows straightforwardly from \cite[Theorem
	5.3.2]{deHaan2007} and is useful in the proof of Theorem
	\ref{theorem:monte}.
	\begin{lemma}
		\label{lemma:moments}
		Let $X_1, \ldots, X_n$ be i.i.d.\ copies of an almost surely
		non-negative $X\in RV_\gamma$ and let $0 < p < 1/\gamma$. Then 
		\begin{equation*}
			\lim_{n\to\infty}
			\mathbb{E}\left(\frac{\max_{i\in\{1, \ldots, n\}} X_i}{U_X(n)}
			\right)^p = \mathbb{E}\left(G^p\right),
		\end{equation*}
		where $G$ follows the Fr\'echet distribution.
	\end{lemma}

	\begin{proof}[Proof of Theorem \ref{theorem:monte}]
		Suppose that the assumptions of Theorem \ref{theorem:monte} hold.
		Similarly as in the proof of Theorem \ref{theorem:continuity}, we start
		by showing that
		\begin{equation}\label{eq:error_monte}
			E_{n,n}
			= O_{\mathbb{P}}\left(\left(\frac{n}{m}\right)^{\frac{1}{2p}}
			U_V(n)\right)
			= O_{\mathbb{P}}\left(U_{\|Y\|_p}(n/k)e'_n\right).
		\end{equation}
		This, together with Theorem \ref{theorem:approx}, gives Parts
		\ref{part1-monte}-\ref{part3-monte}. First observe that
		\eqref{eq:error_monte} follows if we can show that, for any $z>0$, we
		have
		\begin{equation} \label{eq:monte-error-bound2}
			\mathbb{P}\left(\frac{E_{n,n}}{U_{\|Y\|_p}(n/k)} > z\right)
			\leq \left(1 + \varepsilon\right)\left(\frac{1}{z}\right)^{2p}
			\mathbb{E}\left(G^{2p}\right) (e_n')^{2p}.
		\end{equation}
		Set $\delta_n = zU_{\|Y\|_p}(n/k)$. By writing the norm as conditional
		expectation (cf. Remark \ref{remark:monte}) and using the conditional
		Markov inequality together with $|x^{1/p} - y^{1/p}| \leq |x - y|^{1/p}$
		yields
		\begin{equation*}
			\begin{split}
				\mathbb{P}\left(E_i > \delta_n | Y_i\right)
				&= \mathbb{P}\left(\left|\left(\frac{1}{m}\sum_{j=1}^m
				|Y_i(T_j)|^p\right)^{1/p}
				- \left(\mathbb{E}\left(|Y_i(T)|^p | Y_i\right)\right)^{1/p}
				\right| > \delta_n \Bigg| Y_i\right) \\
				&\leq \mathbb{P}\left(\left(\frac{1}{m}\sum_{j=1}^m |Y_i(T_j)|^p
				- \mathbb{E}\left(|Y_i(T)|^p | Y_i\right)\right)^2
				> \delta_n^{2p} \Bigg| Y_i\right) \\
				&\leq \frac{\mathbb{E}\left(\left(\frac{1}{m}\sum_{j=1}^m
				|Y_i(T_j)|^p
				- \mathbb{E}\left(|Y_i(T)|^p | Y_i\right)\right)^2\Bigg|
				Y_i\right)}{\delta_n^{2p}} \\
				&= \frac{\mathbb{E}\left( \var\left(\sum_{j=1}^m |Y_i(T_j)|^p
				\right) \Bigg| Y_i\right)}{m^2\delta_n^{2p}}
				= \frac{\var\left(|Y_i(T)|^p | Y_i\right)}{m\delta_n^{2p}}
				\leq \frac{V_i^{2p}}{m\delta_n^{2p}}.
			\end{split}
		\end{equation*}
		It follows that 
		\begin{equation*}
			\begin{split}
	 			\mathbb{P}\left(E_{n,n} > \delta_n
				\big| Y_1, \dots, Y_n\right)
				&= \mathbb{P}\left(\bigcup_{i=1}^n \left\{E_i > \delta_n\right\}
				\Big| Y_1, \ldots, Y_n\right) \\
				&\leq \sum_{i=1}^n \mathbb{P}\left(E_i > \delta_n\big| Y_i
				\right)
				\leq \sum_{i=1}^n\frac{V_i^{2p}}{m\delta_n^{2p}}
				\leq \frac{n V_{n,n}^{2p}}{m\delta_n^{2p}}.
			\end{split}
		\end{equation*}
		Hence
		\begin{equation*} 
			\begin{split}
				\mathbb{P}\left(E_{n,n} > \delta_n\right)&=\mathbb{E}\left(
				\mathbb{P}\left(E_{n,n} > \delta_n | Y_1, \dots, Y_n\right)\right)
				\leq \frac{n}{m\delta_n^{2p}}
				\mathbb{E}\left(V_{n,n}^{2p}\right) \\
				&= \frac{n}{m} \left(\frac{U_V(n)}
				{U_{\|Y\|_p}(n/k) z}\right)^{2p}
				\mathbb{E}\left(\frac{V_{n,n}}
				{U_V(n)}\right)^{2p} \\
				&= \left(\frac{1}{z}\right)^{2p}
		 		(e_n')^{2p} \mathbb{E}\left(\frac{V_{n,n}}
				{U_V(n)}\right)^{2p}
			\end{split}
		\end{equation*}
		and, as  $\mathbb{E}\left(\frac{V_{n,n}} {U_V(n)}\right)^{2p} \to
		\mathbb{E}\left(G^{2p}\right)$ by Lemma \ref{lemma:moments}, this now
		gives \eqref{eq:monte-error-bound2}.

		It remains to prove Part \ref{part4-monte}. By splitting
		$\mathbb{P}\left(\left|\hat\gamma_n  - \tilde\gamma_n\right| > x\right)$
		into two parts exactly as in the proof of Theorem
		\ref{theorem:continuity} and by applying Proposition
		\ref{prop:min-concentration}, we obtain
		\begin{equation*}
			\begin{split}
				\mathbb{P}\left(\left|\hat\gamma_n  - \tilde\gamma_n\right|
				> x\right)
				&\leq \left(1 + \varepsilon\right)
				\left(\left(a_{y,z}\right)^{-1/(2\gamma)}
				b_{\gamma, \varepsilon, k}\left(a_{y,z}\right)
				\exp\left(a_{y,z}^{-1/\gamma} \frac{k}{k+1} \right)\right) \\
				&+ \mathbb{P}\left(\frac{E_{n,n}}{U_{\|Y\|_p}(n/k)} > z\right)
				+ \mathbb{P}\left(\frac{E_{n,n}}{U_{\|Y\|_p}(n/k)}
				> \frac{x}{2(1+y)}\right).
			\end{split}
		\end{equation*}
		Applying \eqref{eq:monte-error-bound2} two times completes the proof.
	\end{proof}

	\section{Details related to Example \ref{example:riemann} and Example
	\ref{example:monte}}
	\label{appendix:sec:examples}

	We begin by reviewing known results related to sums and products involving
	heavy-tailed random variables. 
	
	\begin{lemma}[\cite{resnick2007}, Proposition 7.5] \label{lemma:breiman} Let
	  	$X$ and $Z$ be nonnegative random variables that are independent of each
	  	other. Assume that $X\in RV_\gamma$ and that for some $\varepsilon > 0$
	  	we have $\mathbb{E}\left(Z^{1/\gamma + \varepsilon}\right) < \infty$.
	  	Then
	  	\begin{equation*}
		  	\lim_{x\to\infty}\frac{\mathbb{P}\left(XZ > x\right)}
		  	{\mathbb{P}\left(X > x\right)} = \mathbb{E}\left(Z^{1/\gamma}
			\right).
	  	\end{equation*}
	\end{lemma}

	\begin{lemma}[\cite{resnick1987}, Lemma 4.24.] \label{lemma:series_rv} Let
		$X_1, X_2, \ldots$ be i.i.d.\ copies of $X$ of a random variable $X$
		satisfying  $\left|X\right|\in RV_\gamma$  and
		\begin{equation*}
			\lim_{x\to\infty}\frac{\mathbb{P}\left(X > x\right)}
			{\mathbb{P}\left(\left|X\right| > x\right)} = p, \quad
			\lim_{x\to\infty}\frac{\mathbb{P}\left(X \leq -x\right)}
			{\mathbb{P}\left(\left|X\right| > x\right)} = q,\quad
			0\leq p\leq 1, \quad p + q = 1.
		\end{equation*}
		Let $\varphi_1, \varphi_2, \ldots$ be a real-valued sequence such that
		$\sum_{j=1}^\infty \left|\varphi_j\right|^\lambda < \infty$ for some
		$ 0 < \lambda < \min\{1/\gamma, 1\}$. Then
		\begin{equation*}
			\lim_{x\to\infty}
			\frac{\mathbb{P}\left(\sum_{j=1}^\infty
			\left|\varphi_j\right|\left|X_j\right| > x\right)}
			{\mathbb{P}\left(\left|X\right| > x\right)}
			= \sum_{j=1}^\infty \left|\varphi_j\right|^{1/\gamma}.
		\end{equation*}
	\end{lemma}

	The next result -- used in Example \ref{example:riemann} -- is useful in
	analyzing H\"older continuity of sample paths. The following Lemma is a
	simplified version of Garsia-Rodemich-Rumsey inequality~\cite[Lemma
	1.1]{garsia1970}.
  	\begin{lemma}\label{lemma:garsia} Let $\theta,\alpha > 0$ such that
		$\alpha\theta > 1$ and let $f$ be a continuous function on $[0, 1]$.
		Then for all $t,s\in [0,1]$,
		\begin{equation*}
	  		|f(t)-f(s)| \leq K_{\alpha, \theta} |t-s|
			^{\alpha - \frac{1}{\theta}}
	  		\left(\int_0^1 \int_0^1 \frac{|f(u)-f(v)|^\theta}
	  		{|u-v|^{\alpha \theta+1}}\,\mathrm{d}u \,\mathrm{d}v\right)
	  		^\frac{1}{\theta},
		\end{equation*}
		where $K_{\alpha, \theta} > 0$ is a constant depending on $\alpha$ and
		$\theta$.
  	\end{lemma}

  	\paragraph*{Details related to Example \ref{example:riemann}.} Suppose that
  	the assumptions given in Example \ref{example:riemann} hold. Our aim is to
  	clarify \eqref{eq:criterion_example}. For this, we begin by showing  that
    \begin{equation} \label{eq:continuity_z}
	  	\left|Z_j(t) - Z_j(s)\right| \leq M_j \left|t - s\right|^\eta
	  	\quad \forall t,s \in [0, 1],
	  	\quad\forall \eta\in(0, \beta/\kappa),
  	\end{equation}
	where $M_j$ are random variables such that $\mathbb{E} M_j^{\kappa} <
	\infty$ and $\mathbb{E} M_j^{1/\gamma'} > 0$. By \eqref{eq:moments} and
	Kolmogorov continuity criterion, $Z$ is (or has a version that is) H\"older
	continuous of any order $\eta\in (0,\beta/\kappa)$. Hence it suffices to
	show the claimed integrability for the H\"older coefficient $M_j$. For this,
	let $\delta\in(0,1)$, $\alpha = ((1 - \delta)\beta + 1)/\kappa$, and $\eta
	= \alpha - 1/\kappa$. Since now $\alpha \kappa>1$, Lemma \ref{lemma:garsia}
	gives
	\begin{equation*}
		|Z(t)-Z(s)|\leq |t-s|^{\eta}K_{\alpha, \kappa} \left(\int_0^1\int_0^1
		\frac{|Z(u) - Z(v)|^{\kappa}}{|u - v|
		^{\alpha\kappa + 1}} \,\mathrm{d}u \,\mathrm{d}v\right)
		^{1/\kappa}.
	\end{equation*}
	Hence, by \eqref{eq:moments} and Tonelli's theorem,
	\begin{equation*}
		M_j = \sup_{0\leq t < s \leq 1} \frac{|Z_j(t) - Z_j(s)|}{|t-s|^{\eta}}
	\end{equation*}
	satisfies
  	\begin{equation*}
	  	\begin{split}
	  		\mathbb{E} M_j^{\kappa}
	  		&\leq K_{\alpha, \kappa}^\kappa
	  		\int_0^1\int_0^1\frac{\mathbb{E}\left[|Z_j(u) - Z_j(v)|
	  		^{\kappa}\right]}{|u - v|^{\alpha\kappa + 1}}
	  		\,\mathrm{d}u \,\mathrm{d}v \\
	  		&\leq K\cdot K_{\alpha, \kappa}^\kappa
	  		\int_0^1\int_0^1\frac{|u - v|^{1 + \beta}}{|u - v|
	  		^{\alpha\kappa + 1}} \,\mathrm{d}u \,\mathrm{d}v.
	  	\end{split} 
  	\end{equation*}
  	This is finite as $(1 + \beta) - (\alpha\kappa + 1) = \delta \beta - 1 >
  	-1$. Since $\delta \in (0,1)$ is arbitrary and since now $M_j>0$ almost
  	surely by our assumptions, we also get  $\mathbb{E} M_j^{1/\gamma'} > 0$ and
  	\eqref{eq:continuity_z} follows. 
  
  	Since in Example \ref{example:riemann} $\{\mathcal{R}_j\}_{j\in\mathbb{N}}$
  	and $\{Z_j\}_{j\in\mathbb{N}}$ are mutually independent, so are
  	$\{\mathcal{R}_j\}_{j\in\mathbb{N}}$ and $\{M_j\}_{j\in\mathbb{N}}$. It then
  	follows that 
  	\begin{equation*}
	  	\begin{split}
		  	&\left|Y(t) - Y(s)\right|
		  	\leq \sum_{j=1}^\infty \left|\varphi_j \mathcal{R}_j\right|
		  	\left|Z_j(t) - Z_j(s)\right| \leq \left|t - s\right|^\eta
			\sum_{j=1}^\infty \left|\varphi_j\right|
		  	\underbrace{M_j \mathcal{R}_j}
		  	_{=: \tilde{\mathcal{R}}_j},
	  	\end{split}
  	\end{equation*}
	where, thanks to Lemma \ref{lemma:breiman},
	$\mathbb{P}\left(\tilde{\mathcal{R}}_j > t\right) \sim
	\mathbb{E}\left(M_j^{1/\gamma'}\right) \mathbb{P}\left(\mathcal{R}_j >
	t\right)$ as $ t\to\infty$. Hence $\tilde{\mathcal{R}}_j \in RV_{\gamma'}$.
	To conclude, Lemma \ref{lemma:series_rv} yields
  	\begin{equation*}
	  	\mathbb{P}\left(\sum_{j=1}^\infty \left|\phi_j\right|
	  	\tilde{\mathcal{R}}_j > t\right)
	  	\sim \mathbb{P}(\tilde{\mathcal{R}}_1 > t)\sum_{j=0}^\infty
	  	\left|\phi_j\right|^{1/\gamma'}, \quad t\to\infty.
	\end{equation*}
 	Thus \eqref{eq:criterion_example} holds with $V:= \sum_{j=1}^\infty
 	\left|\phi_j\right|\tilde{\mathcal{R}}_j \in RV_{\gamma'}$.

  	\paragraph*{Details related to Example \ref{example:monte}.} Suppose that
  	the assumptions given in Example \ref{example:monte} hold. Our aim is to
  	show that $\|Y\|_{2p} \in RV_{\gamma'}$. This follows straightforwardly
  	since Lemma \ref{lemma:breiman} gives
  	\begin{equation*}
  		\mathbb{P}\left(\left\|Y\right\|_{2p} > x\right)
  		= \mathbb{P}\left(\mathcal{R}\left\|Z\right\|_{2p} > x\right)
  		\sim \mathbb{E}\left(\left\|Z\right\|_{2p}^{1/\gamma'}\right)
  		\mathbb{P}\left(\mathcal{R} > x\right), \quad x\to\infty,
  	\end{equation*}
  	where $\mathbb{E}\left(\left\|Z\right\|_{2p}^{1/\gamma'}\right) > 0$ as
  	$\mathbb{P}\left(Z \equiv 0\right) < 1$.
\end{appendices}
\end{document}